\documentclass[12pt]{article}
\usepackage{amscd}
\usepackage{latexsym}
\usepackage{mathrsfs}
\usepackage{pst-all,multido,ifthen}

\usepackage[francais,english]{babel}
\usepackage[intlimits]{amsmath}
\usepackage{amssymb}
\usepackage{latexsym}
\usepackage{mathrsfs}
\usepackage[all]{xy}

\RequirePackage{ifthen}
\RequirePackage{verbatim}

\newcounter{segno}[section]
\renewcommand{\thesegno}{\thesection.\arabic{segno}}

\newcounter{thmno}[section]
\renewcommand{\thethmno}{\thesection.\arabic{thmno}}

\newskip\segskipamount
\newskip\procskipamount
\newskip\interskipamount
\newskip\exskipamount
\newskip\refskipamount
\newcommand{\segskip}{\vskip\segskipamount}
\newcommand{\procskip}{\vskip\procskipamount}
\newcommand{\interskip}{\vskip\interskipamount}
\newcommand{\exskip}{\vskip\exskipamount}
\newcommand{\refskip}{\vskip\refskipamount}

\newcommand{\segbreak}{\par
   \ifdim\lastskip&lt;\segskipamount\removelastskip
   \penalty-200
   \segskip\fi}

\newcommand{\procbreak}{\par
   \ifdim\lastskip&lt;\procskipamount\removelastskip
   \penalty-100
   \procskip\fi}

\newcommand{\exbreak}{\par
   \ifdim\lastskip&lt;\exskipamount\removelastskip
   \penalty-100
   \exskip\fi}

\newcommand{\titlebreak}{\par%
\ifdim\lastskip&lt;\interskipamount\removelastskip%
\penalty10000%
\interskip\fi%
\noindent}%

\newcommand{\interbreak}{\par%
\ifdim\lastskip&lt;\interskipamount\removelastskip%
\penalty-100%
\interskip\fi%
\noindent\ignorespaces}%

\newcommand{\introbreak}{\par
   \ifdim\lastskip&lt;\segskipamount\removelastskip
   \nobreak
   \segskip\fi}

\newcommand{\refbreak}{\par
   \ifdim\lastskip&lt;\refskipamount\removelastskip
   \penalty-100
   \refskip\fi}

\numberwithin{equation}{segno}

\newcommand{\marginrule}{\marginpar{\rule[-10.5mm]{1mm}{10mm}}}

{\par}%

{\par}%

  {\par\interskip}%

  {\par\interskip}%

  {\par\segskip}%

\newenvironment{proof}{%
\noindent%
{\em Proof.\ }}%
{\hspace*{\fill}$\Box$\\[0.5\topskip]\@afterindentfalse\noindent}

  {\par}%

\newcounter{listcounter}
\newcounter{deflistcounter}

\newcounter{exlistcounter}

\newskip{\itemsepamount}
\newskip{\topsepamount}
\newenvironment{definitionlist}{%
  \begin{list}
    {\upshape (\alph{deflistcounter})}
    {\setlength{\leftmargin}{18pt}
     \setlength{\rightmargin}{0pt}
     \setlength{\itemindent}{0pt}
     \setlength{\labelsep}{5pt}
     \setlength{\labelwidth}{13pt}
     \setlength{\listparindent}{\parindent}
     \setlength{\parsep}{0pt}
     \setlength{\itemsep}{\itemsepamount}
     \setlength{\topsep}{\topsepamount}
     \usecounter{deflistcounter}}}
  {\end{list}}

% ...

\renewcommand{\hbar}{\bar{h}}

\DeclareMathOperator{\Lie}{Lie}
\DeclareMathOperator{\re}{red} 
\DeclareMathOperator{\End}{End} 
\DeclareMathOperator{\Hom}{Hom} 
\DeclareMathOperator{\Spec}{Spec}

\DeclareMathOperator{\Ker}{Ker}

\DeclareMathOperator{\GL}{\pmb{GL}} 
 
\DeclareMathOperator{\Image}{Im}

\newcommand\addots{\mathinner{\mkern1mu\raise0pt\vbox{\kern7pt\hbox{.}}\mkern2mu\raise3pt\hbox{.}\mkern2mu\raise6pt\hbox{.}\mkern1mu}}

\newcommand\varto[1]{\mathrel{\hbox to #1pt{\rightarrowfill}}}

\usepackage{amsmath}
\usepackage{amssymb}
\usepackage[all]{xy}
\usepackage{babel}
\newtheorem{thm}{Theorem}
\newtheorem{df}{Def\/inition}
\newtheorem{prop}{Proposition}
\newtheorem{cor}{Corollary}
\newtheorem{lemma}{Lemma}

\newtheorem{fact}{Fact}

\def\endproof{$\hfill \square$}
 
\begin{document}
\title{Dimensions of group schemes of automorphisms of truncated Barsotti--Tate groups}
\author{Ofer Gabber and Adrian Vasiu}
\maketitle

\centerline{Final version, to appear in IMRN}

\bigskip
\noindent
{\bf Abstract.} Let $D$ be a $p$-divisible group over an algebraically closed field $k$ of characteristic $p>0$. Let $n_D$ be the smallest non-negative integer such that $D$ is determined by $D[p^{n_D}]$ within the class of
$p$-divisible groups over $k$ of the same codimension $c$ and dimension $d$ as $D$. We study $n_D$, lifts of $D[p^m]$ to truncated
Barsotti--Tate groups of level $m+1$ over $k$, and the numbers
$\gamma_D(i):=\dim(\pmb{Aut}(D[p^i]))$. We show that  $n_D\le cd$, $(\gamma_D(i+1)-\gamma_D(i))_{i\in\Bbb N}$ is a decreasing sequence in $\Bbb N$, for $cd>0$ we have $\gamma_D(1)<\gamma_D(2)<\cdots<\gamma_D(n_D)$, and for $m\in\{1,\ldots,n_D-1\}$ there exists an infinite set of truncated
Barsotti--Tate groups of level $m+1$ which are pairwise non-isomorphic and lift $D[p^m]$. Different generalizations to $p$-divisible groups with a smooth integral group scheme in the crystalline context are also proved.



\section{Introduction}\label{Sect1}
Let $p\in\Bbb N$ be a prime. Let $k$ be an algebraically closed field of
characteristic $p$. Let $c,d\in\Bbb N$ be such that $r:=c+d>0$. Let $m\in\Bbb
N^{\ast}$. By a $BT_m$ over a scheme $S$ we mean a {\it truncated Barsotti--Tate
group} of level $m$ over $S$ which has codimension $c$ and dimension $d$. In this
paper we put forward {\it principles} that govern: (i) the classification of
$BT_m$'s over $k$ and (ii) the specializations of $BT_m$'s over an algebraic closure
$k_1$ of $k((x))$ to $BT_m$'s over $k$. 

Let $D$ be a {\it $p$-divisible group} over $k$ of codimension $c$ and dimension
$d$. Its height is $r$. Let $n_D$ be the {\it $i$-number} of $D$, i.e. the smallest non-negative integer for which the following statement holds: if $C$ is a
$p$-divisible group over $k$ of codimension $c$ and dimension $d$ and such that
$C[p^{n_D}]$ is isomorphic to $D[p^{n_D}]$, then $C$ is isomorphic to $D$. We have
$n_D=0$ if and only if $cd=0$. For the existence of $n_D$ we refer to [Ma, Ch. III,
Sect. 3], [Tr1, Thm. 3], [Tr2, Thm. 1], [Va1, Cor. 1.3], or [Oo2, Cor. 1.7]. One has
a first estimate $n_D\le cd+1$, cf. [Tr1, Thm. 3]. {\it Traverso's isomorphism
conjecture} predicts that $n_D\le\min\{c,d\}$, cf. [Tr3, Sect. 40, Conj. 4]. 

Thus to classify all $BT_m$'s over $k$ and to decide which $BT_m$'s over $k_1$
specialize to which $BT_m$'s over $k$, one is led to consider four basic problems.

\medskip
{\bf Problem 1.} Determine the set $N_{c,d}$ of possible values of $n_D$.  

\smallskip
{\bf Problem 2.} For $1\le m\le\max(N_{c,d})$ classify all $BT_m$'s over $k$.

\smallskip
{\bf Problem 3.} Find invariants that go up under specializations of $BT_m$'s.

\smallskip
{\bf Problem 4.} Find basic principles that govern {\it group actions} whose {\it orbits} parametrize isomorphism classes of $BT_m$'s over algebraically closed fields. 

\medskip
Significant progress toward the solution of Problem 1 was made in [Va3]. Thus in
this paper we will  present progress toward the solution of Problems 2 to 4
and will provide basic tools for future extra progress toward the solution of Problem 1 (see already \cite{LNV}). We begin by addressing Problem 3.

\subsection{Group schemes of automorphisms}\label{SSect1.1}
Let $B$ be a $BT_m$ over $k$. Let $\pmb{Aut}(B)$ be the group scheme over $k$ of
automorphisms of $B$ and let $\gamma_B:=\dim(\pmb{Aut}(B))$. The importance of the
numbers $\gamma_B$ stems out from the following three main facts:

\medskip
$\bullet$ They are codimensions of versal {\it level $m$ strata}, cf.  [Va2, Thm. 1.2 (c)].

$\bullet$ They can compute the $i$-numbers $n_D$ (see Theorem \ref{Gamma} below). 

$\bullet$  They are a main source of invariants that go up under specializations. 

\medskip\noindent
More precisely, we have the following elementary fact:

\medskip
{\bf (a)} Let $\mathcal R$ be a local domain which is a $k$-subalgebra of $k_1$ and
has residue field $k$. If $\mathcal B$ is a $BT_m$ over $\mathcal R$ that lifts $B$
(i.e. $\mathcal B_k=B$), then $\gamma_{\mathcal B_{k_1}}\le \gamma_B$.

\medskip
Besides $\gamma_B$'s, we know only two other invariants that go up under
specializations of $BT_m$'s. Unfortunately, they are defined only for some specific
values of $m$. To describe them,  let $j:=\lceil {{cd}\over r}\rceil$. It is known
that $j$ is the smallest non-negative integer such that the Newton polygon of each
$p$-divisible group $C$ over $k$ of codimension $c$ and dimension $d$, is uniquely
determined by the isomorphism class of $C[p^j]$ (cf. [NV2, Thm. 1.2]). Thus if $m\ge
j$, we can speak about the Newton polygon $\nu_B$ of $B$: it is the Newton
polygon of every $p$-divisible group over $k$ whose associated $BT_m$ is isomorphic
to $B$. 

\medskip
{\bf (b)} We assume that $m\ge j$. Let $\mathcal R$ and $\mathcal B$ be as in (a).
Then a classical theorem of Grothendieck implies that $\nu_{\mathcal
B_{k_1}}$ is below $\nu_B$. To check this we can assume that $\mathcal R$ is
noetherian. There exists a $p$-divisible group $\widehat{\mathfrak{D}}$ over the
completion $\widehat{\mathcal R}$ of $\mathcal R$ such that
$\widehat{\mathfrak{D}}[p^m]=\mathcal B_{\widehat{\mathcal R}}$, cf. [Il, Thm. 4.4
e)].  Thus the statement on Newton polygons follows from [Ka, Cor. 2.3.2] applied to
the Dieudonn\'e $F$-crystal of $\widehat{\mathfrak{D}}$.

\medskip
{\bf (c)} We assume that $m=1$. Then the set of isomorphism classes of $BT_1$'s over either
$k$ or $k_1$ is in bijection to the finite quotient set $S_r/(S_c\times S_d)$.
Thus one can define naturally an invariant $w_B\in S_r/(S_c\times S_d)$. Moreover
there exists a well understood relation of order $\leqslant $ on $S_r/(S_c\times S_d)$
such that if $\mathcal R$ and $\mathcal B$ are as in (a), then we have $w_{\mathcal
B_{k_1}}\leqslant w_B$. See [We] and [PWZ] for all these facts (see also [Kr], [Oo1], [Mo], [MW], and [Va4]). We emphasize that the invariant $w_B$ is a refinement of the classical invariants of $B$ (like the {\it $a$-number} of $B$, the {\it $p$-rank} of $B$, the $p$-rank of the {\it Cartier dual} $B^{\text{t}}$ of $B$, etc.).

\medskip
We recall that the $a$-number $a_D$ of $D[p]$ (or of $D$) is the greatest non-negative integer such that $\pmb{\alpha}_p^{a_D}$ is a subgroup scheme of $D[p]$. For $i\in\Bbb N$ let $\gamma_D(i):=\dim(\pmb{Aut}(D[p^i]))$; thus $\gamma_D(0)=0$.
We call $(\gamma_D(i))_{i\in\Bbb N}$ the {\it centralizing sequence} of $D$ and we
call $s_D:=\gamma_D(n_D)$ the {\it specializing height} of $D$.

\begin{thm}\label{Gamma}
The centralizing sequence of $D$ has the following three basic properties:

\medskip
{\bf (a)} For each $l\in\Bbb N$, the sequence
$(\gamma_D(i+l)-\gamma_D(i))_{i\in\Bbb N}$ is a decreasing sequence in $\Bbb N$.

\smallskip
{\bf (b)} If $cd>0$, then we have $a_D^2\le \gamma_D(1)<\gamma_D(2)<\cdots <\gamma_D(n_D)$.

\smallskip
{\bf (c)} For $i\ge n_D$ we have $\gamma_D(i)=\gamma_D(n_D)\le cd$.
\end{thm}

\medskip
The proof of Theorem \ref{Gamma} (a) relies on properties of the group scheme $\pmb{End}(D[p^i])$ over $k$ of
endomorphisms of $D[p^i]$ presented in Subsection \ref{SSect2.1} (see Subsection \ref{SSect2.2}). The
proof of Theorem \ref{Gamma} (b) relies on Theorem \ref{Gamma} (a) and on properties
of the group actions of Section 3 (see Subsection \ref{PGamma}). If $n_D>0$ (i.e. if
$cd>0$), then Theorem \ref{Gamma} (c) was proved in [Va2, Thm. 1.2 (d) and Cor.
1.2.1]. The case $n_D=0$ of Theorem \ref{Gamma} (c) is trivial.

The identity component (resp. the reduced group) of an affine group scheme $\sharp$ over $k$ will be denoted as $\sharp^0$ (resp. $\sharp_{\re}$). A connected smooth affine group over $k$ is called {\it unipotent} if $\Bbb G_m$
is not a subgroup of it. 

The following three consequences of Theorem \ref{Gamma} and its proof are proved in Subsections \ref{Pcor1} to \ref{Pcor3}  (respectively).

\begin{cor}\label{cor1}
We have $n_D\le cd$. If $D$ is not ordinary, then we have $n_D\le cd+1-a_D^2\le cd$. 
\end{cor}

\begin{cor}\label{cor2} The following three properties hold:

\medskip
{\bf (a)} Let $i,l\in\Bbb N^{\ast}$ with $i<l$. Then the image of the restriction homomorphism $\pmb{Aut}(D[p^l])\to \pmb{Aut}(D[p^i])$ (or of $\pmb{End}(D[p^l])\to \pmb{End}(D[p^i])$) is finite if and only if $l-i\ge n_D$. 

\smallskip
{\bf (b)} We assume that $D$ is not ordinary. Let $n$ be the greatest positive
integer such that the image of the restriction homomorphism
$\pmb{Aut}(D[p^n])\to \pmb{Aut}(D[p])$ (or of $\pmb{End}(D[p^n])\to \pmb{End}(D[p])$) has positive dimension. Then we have $n=n_D$.

\smallskip
{\bf (c)} For all $s\in\Bbb N^{\ast}$ we have $n_{D^s}=n_D$.
\end{cor}

\begin{cor}\label{cor3} Let $i,l\in\Bbb N^{\ast}$ with $i<l$. Then the
following five properties hold:

{\bf (a)} The restriction homomorphism $\pmb{Aut}(D[p^l])_{\re}^0\to \pmb{Aut}(D[p^i])_{\re}^0$
is trivial if and only if $l-i\ge n_D$.

\smallskip
{\bf (b)} The connected smooth affine group scheme $\pmb{Aut}(D[p^i])_{\re}^0$ is unipotent. 

\smallskip
{\bf (c)} If $i\ge n_D+1$, then the exponent of the group $\pmb{Aut}(D[p^i])^0(k)$ is $p^{n_D}$.

\smallskip
{\bf (d)} If $i\ge 2n_D$, then $\pmb{Aut}(D[p^i])_{\re}^0$ is commutative. 

\smallskip
{\bf (e)} If $i\ge n_D+1$, then there exist canonical bijections
$\pmb{Aut}(D[p^l])^0(k)\to \pmb{Aut}(D[p^i])^0(k)$. If moreover $i\ge 2n_D$, then these bijections
are group isomorphisms.
\end{cor}

Corollary \ref{cor1} refines [Tr1, Thm. 3] (if $D$ is not ordinary, then $a_D$ could be any integer in the interval $[1,\min\{c,d\}]$). See Subsection \ref{Rcor1} for a
refinement of Corollary \ref{cor1} in terms of Newton polygon slopes of $D$. Based on the above results and [Va3], in \cite{LNV} a corrected, refined, and optimal version of Traverso's isomorphism conjecture is proved. 

\subsection{Applications of group actions}\label{SSect1.2}
In Section 3 we introduce certain group actions over $k$ whose sets of orbits
parametrize isomorphism classes of (truncations of) $p$-divisible groups over $k$ of
codimension $c$ and dimension $d$. These group actions {\it encode all the data one
could possible think of}  in connection to (truncations of) $p$-divisible groups  over $k$  of codimension $c$ and dimension $d$ (e.g., encode their $i$-numbers,
their group schemes of automorphisms, their specializations, all possible
stratifications associated to them, their local deformations, etc.). As in [Tr3,
Sects. 26 to 39], [Va2, Sect. 2], and [Va4, Subsect. 4.2.2 and Sect. 5.1], these
group actions involve (projective limits of) connected smooth affine groups
over $k$. The idea of using group actions to study isomorphism classes of
(truncations of) $p$-divisible groups over $k$ goes back to Manin (see [Ma, Ch. III,
Sect. 3]). But Manin used disjoint unions (indexed by finite sets of isomorphism
classes of special $p$-divisible groups over $k$) of discrete actions (i.e. of
actions involving finite groups). In this paper, we only use group actions in order
to prove the two principles below and to apply the first principle to the proofs of
Theorem \ref{Gamma} (b) and (c) and of Theorem \ref{INF} below.

We call two group schemes over $k$ {\it radicially
isogenous} if their perfections are isomorphic group schemes over $k$ (to be
compared with [Se, Subsect. 1.2, Def. 1], where the terminology {\it equivalent} is
used).  The next two theorems are proved in Subsections  \ref{PTUP} and \ref{PCP} (respectively)
and form progress towards the solution of Problems 2 and 4. 

\begin{thm}{\bf (Unipotent Principle)}\label{UP}
Let $c,d, m$, and $D$ be as above. Then the set of isomorphism classes of lifts of
$D[p^m]$ to $BT_{m+1}$'s over $k$ is in natural bijection to the set of orbits of a
suitable action $\mathcal T_{m+1}\times_k\mathcal V_{m+1}\to\mathcal V_{m+1}$ over
$k$ which has the following three properties:

\medskip
{\bf (i)} the group scheme $\mathcal T_{m+1}$ over $k$ is smooth, affine, and its
identity component is a unipotent group of dimension $r^2+\gamma_D(m)$;

\smallskip
{\bf (ii)} the $k$-scheme $\mathcal V_{m+1}$ is an affine space of dimension $r^2$;

\smallskip
{\bf (iii)} we have a short exact sequence $1\to\Bbb G_a^{r^2}\to \mathcal
T_{m+1}\to\mathcal C_m\to 1$, where $\mathcal C_m$ is a smooth, affine group over
$k$ that is radicially isogenous to $\pmb{Aut}(D[p^m])_{\re}$.
\end{thm}

\begin{thm}{\bf (Centralizing Principle)}\label{CP} 
Let $D$ be as above. Let $\mathcal R$ be the normalization of $k[[x]]$ in an
algebraic closure $k_1$ of $k((x))$. Let $\mathfrak{D}$ be a $p$-divisible group
over $\Spec \mathcal R$ whose fibre over $k$ is $D$. We assume that there exists a $p$-divisible group $D^\prime$ over $k$ such that $D^\prime_{k_1}$ is isomorphic to $\mathfrak{D}_{k_1}$ (thus we have $\gamma_{D^\prime}(i)\le\gamma_D(i)$ for all $i\in\Bbb N^{\ast}$, cf. fact (a) of Subsection \ref{SSect1.1}). Then the following two properties hold:

\medskip
{\bf (a)} If there exists $l\in\Bbb N^{\ast}$ such that
$\gamma_{D^\prime}(l)=\gamma_D(l)$, then $D^\prime[p^l]$ is isomorphic to $D[p^l]$;
thus for $i\in\{1,\ldots,l\}$ we have $\gamma_{D^\prime}(i)=\gamma_D(i)$.

\smallskip
{\bf (b)} We have $s_{D^\prime}\le s_D$. If $s_{D^\prime}=s_D$, then $D^\prime$ is isomorphic to $D$.
\end{thm}

One could phrase the Centralizing Principle as: {\it centralizing sequences go
quasi-strictly up under specializations}. Theorem \ref{CP} (b) justifies our
terminology specializing height. In Subsection \ref{PINF} we use the strict
inequalities of Theorem \ref{Gamma} (b) and the Unipotent Principle to prove the
following theorem. 

\begin{thm}\label{INF}
Let $B$ be an arbitrary $BT_m$ over $k$. Let $D$ be a $p$-divisible group over $k$
such that $B=D[p^m]$, cf. [Il, Thm. 4.4 e)]. If $m<n_D$, then there exists an
infinite set of $BT_{m+1}$'s over $k$ which are pairwise non-isomorphic and which
lift $B$.
\end{thm}

The case $m=1$ of Theorem \ref{INF} is a stronger form of [Oo3, Thm. B]. 

\subsection{Generalizations}\label{SSect1.3}
In Section \ref{Sect6} we present different generalizations of Theorem \ref{Gamma} to relative contexts provided by quadruples of the form $(M,\phi,\vartheta,G)$, where $(M,\phi,\vartheta)$ is the contravariant Dieudonn\'e module of $D$ and where $G$ is a smooth integral closed subgroup scheme of $\GL_M$ subject to a weakening of the two axioms of [Va2, Subsect. 4.1] (the case $G=\GL_M$ corresponds to the above classical context provided by $D$ itself). Three examples of special interest are presented in Section \ref{Sect7}. They pertain to group schemes of homomorphisms between $BT_m$'s, to group schemes of automorphisms of filtered $BT_m$'s, and to group schemes of automorphisms of $BT_m$'s equipped with (symmetric) principal quasi-polarizations. The motivation for all these generalizations stems out from applications to level $m$ stratifications of special fibres of good integral models of Shimura varieties of Hodge type in unramified mixed characteristic $(0,p)$ (see [Va2, Sect 4]). 

Appendix A contains basic properties of affine group schemes over discrete valuation rings that are needed to present the group actions of Sections 3 and \ref{Sect6}. Appendix B presents properties of exponential and logarithmic maps for matrices in $p$-adic contexts that are required for the mentioned generalizations. 

\section{Group schemes of endomorphisms}\label{Sect2}

Always $c$, $d$, $r$, $m$, $D$, $a_D$, $\pmb{Aut}(D[p^i])$, and $\pmb{End}(D[p^i])$ will be as in Section 1. Let $l,i\in\Bbb N$ be such that $l\ge
i$. We first list three properties of $\pmb{End}(D[p^i])$ (see Subsection \ref{SSect2.1}) and then we use them to prove
Theorem \ref{Gamma} (a) (see Subsection \ref{SSect2.2}). A homomorphism between affine group
schemes over $k$ is called an epimorphism (resp. a monomorphism) if it is
faithfully flat (resp. is a closed embedding). 

\subsection{Basic properties}\label{SSect2.1} 
As the scheme $\pmb{Aut}(D[p^i])$ is an open subscheme of $\pmb{End}(D[p^i])$ and as
all connected components of $\pmb{End}(D[p^i])$ have the same dimension, we get
that:

\medskip
{\bf (i)} We have $\gamma_D(i)=\dim(\pmb{End}(D[p^i]))$.

\medskip
Let $j_{l,i}:D[p^l]\twoheadrightarrow D[p^i]$ be the epimorphism defined by the
multiplication by $p^{l-i}$. Let $i_{i,l}:D[p^i]\hookrightarrow D[p^l]$ be the
natural inclusion monomorphism. 

We consider the monomorphism 
$$\kappa_{i,l}:\pmb{End}(D[p^i])\hookrightarrow \pmb{End}(D[p^l])$$ 
defined on valued points as follows. If $R$ is a commutative $k$-algebra, then
$\kappa_{i,l}(R)$ maps $f_i\in \pmb{End}(D[p^i])(R)$ to $i_{i,l,R}\circ f_i\circ
j_{l,i,R}\in \pmb{End}(D[p^l])(R)$. Let 
$$r_{l,i}:\pmb{End}(D[p^l])\to \pmb{End}(D[p^i])$$ 
be the restriction homomorphism
(associated to the monomorphism $i_{i,l}:D[p^i]\hookrightarrow D[p^l]$). Let
$f_l\in\pmb{End}(D[p^l])(R)$. It is easy to see that there exists $f_i\in
\pmb{End}(D[p^i])(R)$ such that $\kappa_{i,l}(R)(f_i)=f_l$ if and only if
$\text{Ker}(f_l)$ contains $D[p^{l-i}]_R=\text{Ker}(j_{l,i})_R$ (i.e. if and only
if $r_{l,l-i}(R)(f_l)=0\in \pmb{End}(D[p^{l-i}])(R)$). From this we get that: 

\medskip
{\bf (ii)} We have an exact complex
 $$0\to \pmb{End}(D[p^i])\xrightarrow{\kappa_{i,l}}
\pmb{End}(D[p^l])\xrightarrow{r_{l,l-i}} \pmb{End}(D[p^{l-i}]).$$

From properties (i) and (ii), at the level of dimensions we get that:

\medskip
{\bf (iii)} The image $\Image(r_{l,l-i})$ has dimension
$\gamma_D(l)-\gamma_D(i)$.

\subsection{Proof of Theorem \ref{Gamma} (a)}\label{SSect2.2} Let $i,l\in \Bbb N$. As we have
$r_{i+1+l,l}=r_{i+l,l}\circ r_{i+1+l,i+l}$, $\Image(r_{i+1+l,l})$ is a subgroup
scheme of $\Image(r_{i+l,l})$. Thus we have $\dim(\Image(r_{i+l,l}))\ge
\dim(\Image(r_{i+1+l,l}))$. From this and the property \ref{SSect2.1} (iii) we get that the
sequence $(\gamma_D(i+l)-\gamma_D(i))_{i\in\Bbb N}$ is a decreasing sequence in
$\Bbb N$. Therefore Theorem \ref{Gamma} (a) holds.\endproof

\section{Isomorphism classes of $p$-divisible groups as orbits of group actions}\label{Sect3}

We construct group actions over $k$ whose sets of orbits parametrize
isomorphism classes of $p$-divisible groups over $k$ of codimension $c$ and
dimension $d$ (see Subsections \ref{SSect3.1} and \ref{act} that follow the pattern of [Va2]).
Theorem \ref{UP} is proved in Subsection  \ref{PTUP}. Subsection \ref{orbits} recalls basic results
on group actions over $k$; we use them to prove Theorems \ref{Gamma}, \ref{CP}, and \ref{INF} in Subsections \ref{PGamma}, \ref{PCP}, and \ref{PINF} (respectively). 

Let $W(k)$ be the ring of ($p$-typical) Witt vectors with coefficients in $k$. Let $B(k)$ be the
field of fractions of $W(k)$. Let $\sigma$ be the Frobenius automorphism of $k$,
$W(k)$, and $B(k)$. Let $(M,\phi,\vartheta)$ be the contravariant Dieudonn\'e module
of $D$ in the sense of [BBM, Def. 3.1.5], where $\phi$ and $\vartheta$ are the $F$ and $V$ of loc. cit. Thus $M$ is a free $W(k)$-module of rank $r$, $\phi:M\to M$ is
a $\sigma$-linear endomorphism, $\vartheta:M\to
M$ is a $\sigma^{-1}$-linear endomorphism, and we have $\phi\circ\vartheta=\vartheta\circ\phi=p 1_M$. The $\sigma$-linear automorphism $\phi[{1\over p}]:M[{1\over p}]\to M[{1\over p}]$ and its inverse will be denoted simply by $\phi$ and $\phi^{-1}$ (respectively); thus $\phi^{-1}(M):=(\phi[{1\over p}])^{-1}(M)$. We denote also by $\phi$ the $\sigma$-linear
automorphism of $\End(M)[{1\over p}]$ which takes $e\in\End(M)[{1\over p}]$ to
$\phi(e):=\phi\circ e\circ\phi^{-1}\in \End(M)[{1\over p}]$. We refer to [BBM] for
the crystalline contravariant Dieudonn\'e functor $\Bbb D$ defined on the category
of $p$-divisible groups over $k$-schemes. 

We recall that $m\in\Bbb N^{\ast}$. Let $\phi_m,\vartheta_m:M/p^mM\to M/p^mM$ be the
reductions modulo $p^m$ of $\phi,\vartheta:M\to M$. Let $\natural^{(\sigma^m)}$ be the
pull-back (or the tensorization) of $\natural$ with $\sigma^m$. Thus
$M^{(\sigma)}:=W(k)\otimes_{\sigma,W(k)} M$, etc. 

\subsection{Group schemes}\label{SSect3.1} 

In this subsection we introduce several affine group schemes. 

\subsubsection{A dilatation}\label{dil}
Let $\bar F^1$ be the kernel of $\phi_1$; it is a $k$-vector subspace of
$M/pM$ of dimension $d$. Let $\mathcal W$ be the normalizer of $\bar F^1$ in
$\GL_{M/pM}$. If $cd=0$, then $\mathcal W=\GL_{M/pM}$. If $cd>0$, then
$\mathcal W$ is a maximal proper parabolic subgroup of $\GL_{M/pM}$. Let
$\tilde{\mathcal H}$ be the {\it dilatation} of $\GL_M$ centered on the smooth
subgroup $\mathcal W$ of $\GL_{M/pM}$ (see Appendix A). Thus $\tilde{\mathcal H}$ is a smooth affine group scheme over $\Spec W(k)$ equipped with a homomorphism $\tilde{\mathcal P}:\tilde{\mathcal H}\to
\GL_M$ that has the following universal property: a morphism $f:S\to \GL_M$ of flat $\Spec W(k)$-schemes factors
(uniquely) through $\tilde{\mathcal P}$ if and only if the morphism $f_k:S_k\to
\GL_{M/pM}$ factors through $\mathcal W$. The generic fibre of $\tilde{\mathcal P}$ is an isomorphism of $\Spec B(k)$-schemes and thus the relative dimension of $\tilde{\mathcal H}$ over $\Spec W(k)$ is $r^2$. The special fibre $\tilde{\mathcal H}_k$ is connected and we can redefine $\tilde{\mathcal P}$ as follows (cf. Example of Appendix A for details). We can identify $\tilde{\mathcal H}$ with the schematic closure $\GL_{M,\phi^{-1}(M)}$ in $\GL_M\times_{W(k)} \GL_{\phi^{-1}(M)}$ of $\GL_{M[{1\over p}]}$ embedded diagonally in  $\GL_M\times_{W(k)} \GL_{\phi^{-1}(M)}$. Under this identification, $\tilde{\mathcal P}$ gets identified with the first projection homomorphism $\GL_{M,\phi^{-1}(M)}\to \GL_M$. The homomorphism $\tilde{\mathcal H}\to \GL_{\phi^{-1}(M)}$ which is the composite of the identification $\tilde{\mathcal H}=\GL_{M,\phi^{-1}(M)}$ with the second projection homomorphism $\GL_{M,\phi^{-1}(M)}\to \GL_{\phi^{-1}(M)}$, is isomorphic to the dilatation of $\GL_{\phi^{-1}(M)}$ centered on the parabolic subgroup of $\GL_{\phi^{-1}(M)/p\phi^{-1}(M)}$ that normalizes the $k$-vector subspace $M/p\phi^{-1}(M)$ of $\phi^{-1}(M)/p\phi^{-1}(M)$.

For a commutative $W(k)$-algebra $R$ and for $\tilde h\in \tilde{\mathcal
H}(R)$, let $h:=\tilde{\mathcal P}(R)(\tilde h)\in \GL_M(R)$. For
$e\in\End(M)$, $g\in \GL_M(W(k))$, and $\tilde h\in\tilde{\mathcal H}(W(k))$, let $e[m]\in\End(M/p^mM)$, $g[m]\in\GL_M(W_m(k))$, and $\tilde h[m]\in\tilde{\mathcal H}(W_m(k))$ be the reductions modulo $p^m$ (here $W_m(k)=W(k)/p^mW(k)$).

\begin{fact}\label{F1}
Let $g,g_1\in \GL_M(W(k))$. Let $g_2:=gg_1\phi(g^{-1})\in\GL_M(B(k))$ (thus $gg_1\phi=g_2\phi g$). Then we have $g_2\in \GL_M(W(k))$ if and only if there
exists $\tilde h\in \tilde{\mathcal H}(W(k))$ such that $g=h$.
\end{fact}

\noindent
{\bf Proof:} We have $g_2=gg_1\phi(g^{-1})\in \GL_M(W(k))$ if and only if $g$
normalizes $\phi^{-1}(M)$ and thus if and only if $g[1]$ normalizes
$p\phi^{-1}(M)/pM=\bar F^1$. But $g[1]$ normalizes $\bar F^1$ (i.e. we have
$g[1]\in \mathcal W(k)$) if and only if there exists $\tilde h\in\tilde{\mathcal
P}(W(k))$ such that $g=h$, cf. the above universal property.\endproof

\subsubsection{On $\tilde{\mathcal H}^{(\sigma)}$}\label{tildeH}
The group scheme $\tilde{\mathcal H}^{(\sigma)}$ is the dilatation of
$\GL_{M^{(\sigma)}}$ centered on the subgroup scheme $\mathcal
W^{(\sigma)}$ of $\GL_{M^{(\sigma)}/pM^{(\sigma)}}$. One can define
$\tilde{\mathcal H}^{(\sigma)}$ more intrinsically as follows. We identify $\phi$
and $\vartheta$ with $W(k)$-linear maps $\phi:M^{(\sigma)}\to M$ and
$\vartheta:M\to M^{(\sigma)}$. Let $\pmb{Aut}(\phi,\vartheta)$ be the closed subgroup scheme of $\GL_{M^{(\sigma)}}\times_{W(k)} \GL_M$ that represents the functor of automorphisms of the pair of $W(k)$-linear maps $\phi$ and $\vartheta$. If $R$ is a commutative $W(k)$-algebra, then
$\pmb{Aut}(\phi,\vartheta)(R)$ is the group of those pairs $(\tilde a,\tilde
b)\in\GL_{M^{(\sigma)}}(R)\times\GL_M(R)$ for which we have $\tilde
b\circ (1_R\otimes\phi)=(1_R\otimes\phi)\circ \tilde a$ and $\tilde a\circ
(1_R \otimes\vartheta)=(1_R \otimes\vartheta)\circ \tilde b$. 

Under the $W(k)$-linear isomorphism $M\to \phi^{-1}(M)^{(\sigma)}$ induced by $\phi^{-1}$, the group scheme $\pmb{Aut}(\phi,\vartheta)$ becomes isomorphic to the group scheme over $\Spec W(k)$ that represents the functor of automorphisms of the pair of $W(k)$-linear maps $\xymatrixcolsep{1pc}\xymatrix@1{M^{(\sigma)} \ar@<0.5ex>[rr]^{j^{(\sigma)}} & &\phi^{-1}(M)\ar@<0.5ex>[ll]^{p} \\}$$^{(\sigma)}$, where $j:M\hookrightarrow \phi^{-1}(M)$ is the inclusion. Thus $\pmb{Aut}(\phi,\vartheta)$ is isomorphic to the schematic closure $\GL_{M^{(\sigma)},\phi^{-1}(M)^{(\sigma)}}$ in $\GL_{M^{(\sigma)}}\times_{W(k)} \GL_{\phi^{-1}(M)^{(\sigma)}}$ of $\GL_{M^{(\sigma)}[{1\over p}]}$ embedded diagonally in the product $\GL_{M^{(\sigma)}}\times_{W(k)} \GL_{\phi^{-1}(M)^{(\sigma)}}$, cf. Example of Appendix A.
Thus we have natural isomorphisms $\tilde{\mathcal H}^{(\sigma)}\to\GL_{M^{(\sigma)},\phi^{-1}(M)^{(\sigma)}}\to\pmb{Aut}(\phi,\vartheta)$ (cf. Subsubsection \ref{dil} for the first one) and under their composite, the homomorphism $\tilde{\mathcal P}^{(\sigma)}:\tilde{\mathcal H}^{(\sigma)}\to \GL_M^{(\sigma)}=\GL_{M^{(\sigma)}}$ gets identified with the first projection homomorphism $\pmb{Aut}(\phi,\vartheta)\to\GL_{M^{(\sigma)}}$. The first (resp. second) projection homomorphism $\pmb{Aut}(\phi,\vartheta)\to\GL_{M^{(\sigma)}}$ (resp. $\pmb{Aut}(\phi,\vartheta)\to\GL_M$) is the dilatation of
$\GL_{M^{(\sigma)}}$ centered on the parabolic subgroup $\mathcal
W^{(\sigma)}$ of $\GL_{M^{(\sigma)}/pM^{(\sigma)}}$ (resp. is the dilatation of $\GL_M$ centered on the parabolic subgroup of $\GL_{M/pM}$ which is the normalizer of the kernel of $\vartheta_1$).

\subsubsection{The $\Bbb W_m$ functor}\label{FunctorW_m} Let $\text{Aff}_k$ be the category of affine
schemes over $k$. Let $\text{Set}$ be the category of sets. Let $X$ be an affine
scheme over $\Spec W(k)$ of finite type. Let $\Bbb
W_m(X):\text{Aff}_k\to\text{Set}$ be the contravariant functor that associates to an
affine $k$-scheme $\Spec R$ the set $X(W_m(R))$, where $W_m(R)$ is the ring of Witt
vectors of length $m$ with coefficients in $R$. It is well known that this functor
is representable by an affine $k$-scheme of finite type (cf. [Gre, Sect. 4, Cor. 1,
p. 639]) to be also denoted as $\Bbb W_m(X)$. We have $\Bbb W_m(X)(k)=X(W_m(k))$ as well as a functorial identification $\Bbb W_1(X)=X_k$. If $f:Y\to X$ is a closed embedding, then $\Bbb W_m(f):\Bbb W_m(Y)\to\Bbb W_m(X)$ is also a closed embedding (cf. [Gre, Sect. 4, Cor. 2, p. 639]).

We assume that $X$ is smooth over $\Spec W(k)$. If $J$ is an ideal of $R$ of square
$0$, then the ideal $\Ker(W_m(R)\to W_m(R/J))$ is of square $0$ and therefore the reduction map $X(W_m(R))\to X(W_m(R/J))$ is surjective (cf. [BLR, Ch. 2, Sect. 2.2,
Prop. 6]). From this and loc. cit. we get that the affine $k$-scheme $\Bbb W_m(X)$
of finite type is smooth. 

We assume that $X$ is a smooth affine group scheme over $\Spec W(k)$. Then $\Bbb
W_m(X)$ is a smooth affine group over $k$. The length reduction $W(k)$-epimorphisms
$W_{m+1}(R)\twoheadrightarrow W_m(R)$ define naturally a smooth epimorphism
$$\text{Red}_{m+1,X}:\Bbb W_{m+1}(X)\twoheadrightarrow \Bbb W_m(X)$$ 
of affine groups
over $k$ (the smoothness part is checked as in the previous paragraph). The kernel of $\text{Red}_{m+1,X}$ is the vector group over $k$ defined by
the Lie algebra $\Lie(X_k)^{(\sigma^m)}$, and thus it is a unipotent commutative group isomorphic to $\Bbb G_a^{\dim(X_k)}$. Using this and the
identification $\Bbb W_1(X)=X_k$, by induction on $m\in\Bbb N^{\ast}$ we get that: 

\medskip
{\bf (i)} we have $\dim(\Bbb W_m(X))=m\dim(X_k)$;

\smallskip
{\bf (ii)} the group $\Bbb W_m(X)$ is connected if and only if $X_k$ is connected.

\subsection{The group actions $\Bbb T$}\label{act}

Let $\mathcal H_m:=\Bbb W_m(\tilde{\mathcal H})$ and $\mathcal D_m:=\Bbb
W_m(\GL_{M})$. In this paper we will use only occasionally the fact that $\mathcal D_m$ has a natural group structure over $k$. Thus we will view $\mathcal H_m$ as a connected smooth affine group over $k$ of dimension $mr^2$ and (except in few places) we will view $\mathcal D_m$ as a connected smooth affine variety over $k$ of dimension $mr^2$, cf. Subsubsection  \ref{FunctorW_m}. We have an action 
$$\Bbb T_m:\mathcal H_m\times_k \mathcal D_m\to \mathcal D_m$$
defined on $k$-valued points as follows. If $\tilde h\in\tilde{\mathcal H}(W(k))$ and $g\in \GL_{M}(W(k))$, then the product of $\tilde h[m]\in \mathcal
H_m(k)=\tilde{\mathcal H}(W_m(k))$ and $g[m]\in\mathcal D_m(k)=\GL_M(W_m(k))$
is the element (see Subsubsection \ref{dil} for $h$)
$$\Bbb T_m(\tilde h[m],g[m]):=(hg\phi(h^{-1}))[m].$$
This makes sense, cf. Fact \ref{F1}. The fact that $\Bbb T_m$ is indeed a morphism follows easily from [Va2, Formula (1a)] and its natural analogs for arbitrary commutative $k$-algebras. We recall (see [Va2, Lemma 2.2.1]) the following lemma: 

\begin{lemma}\label{L1}
 Let $g_1,g_2\in \GL_M(W(k))$. Then the points $g_1[m],g_2[m]\in\mathcal
D_m(k)$ belong to the same orbit of the action $\Bbb T_m$ if and only if the
Dieudonn\'e modules $(M/p^mM,g_1[m]\phi_m,\vartheta_m g_1[m]^{-1})$ and 
$(M/p^mM,g_2[m]\phi_m,\vartheta_m g_2[m]^{-1})$ are isomorphic.
\end{lemma}

\subsubsection{The limit process}\label{limit} Let $\rho_{m+1}:=\text{Red}_{m+1,\tilde{\mathcal H}}$ and
$\tau_{m+1}:=\text{Red}_{m+1,\GL_M}$ (see Subsubsection  \ref{FunctorW_m} for notations).
Thus $\rho_{m+1}:\mathcal H_{m+1}\twoheadrightarrow\mathcal H_m$ is the natural
reduction epimorphism of affine groups over $k$ and $\tau_{m+1}:\mathcal
D_{m+1}\twoheadrightarrow\mathcal D_m$ is the natural reduction faithfully flat morphism (epimorphism) of
affine (group) varieties over $k$. The action $\Bbb T_m$ is a natural reduction of
the action $\Bbb T_{m+1}$. In other words, the following identity holds:
$$\tau_{m+1}(\Bbb T_{m+1}(\tilde h[m+1],g[m+1]))=\Bbb T_m(\tilde h[m],g[m]).$$
As $\rho_{m+1}$ and $\tau_{m+1}$ are morphisms of affine schemes, the projective
limits 
$$\mathcal H_{\infty}:=\text{proj.}\text{lim.}_{m\to\infty} \mathcal
H_m\;\;\;\;\;\;\text{and}\;\;\;\;\;\;\mathcal D_{\infty}:=\text{proj.}\text{lim.}_{m\to\infty}
\mathcal D_m$$ 
in the category of ringed spaces are also projective
limits in the category of $k$-schemes (cf. [Gro, Sect. 8, Rm. 8.2.14]), $\mathcal H_{\infty}$ is an affine group
scheme over $k$, and $\mathcal D_{\infty}$ is an affine (group) scheme over $k$. We have a limit action
$$\Bbb T_{\infty}:\mathcal H_{\infty}\times_k\mathcal D_{\infty}\to\mathcal
D_{\infty}$$
with the property that for all elements $\tilde h\in\mathcal
H_{\infty}(k)=\tilde{\mathcal H}(W(k))$ and $g\in\mathcal
D_{\infty}(k)=\GL_M(W(k))$, the product $\Bbb T_{\infty}(\tilde h,g)\in\mathcal
D_{\infty}(k)$ is the unique element such that for all $m\in\Bbb N^{\ast}$ we have
$\Bbb T_{\infty}(\tilde h,g)[m]=\Bbb T_m(\tilde h[m],g[m])$; thus we have $\Bbb
T_{\infty}(\tilde h,g)=hg\phi(h^{-1})\in\mathcal D_{\infty}(k)=\GL_M(W(k))$.

There exist several equivalent ways to define orbits of the action $\Bbb T_{\infty}$ (see [GV]). To be short, in what follows by an orbit of $\Bbb T_{\infty}$ we mean a reduced, quasi-compact, locally closed subscheme $\mathcal O^{\prime}_{\infty}$ of $\mathcal D_{\infty}$ which is a natural projective limit $\text{proj.}\text{lim.}_{m\to\infty}\mathcal O_m^\prime$ in the category of ringed spaces of orbits $\mathcal O_m^\prime$ of $\Bbb T_m$ (i.e. classical orbits of $k$-points) under faithfully flat transition morphisms $\mathcal O_{m+1}^\prime\to\mathcal O_m^\prime$ that are induced by $\tau_{m+1}$. As we have surjective transition morphisms between non-empty varieties over $k$, $\mathcal O_{\infty}^\prime$ has $k$-valued points.

\subsubsection{Notations}\label{not}
Let $\mathcal O_m$ be the orbit of $1_M[m]\in\mathcal
D_m(k)$ under the action $\Bbb T_m$; it is a locally closed subscheme of $\mathcal D_m$ which is integral (as $\mathcal H_m$ is connected). Let $\mathcal S_m$ be the subgroup scheme of $\mathcal H_m$ which is the stabilizer
of $1_M[m]$. Let $\mathcal C_m:=\mathcal S_{m,\re}$. We have 
$$\dim(\mathcal S_m)=\dim(\mathcal C_m)=\dim(\mathcal C_m^0)=\dim(\mathcal H_m)-\dim(\mathcal O_m).$$
\indent
Let $\mathcal T_{m+1}$ be the reduced group of the group subscheme
$\rho_{m+1}^{-1}(\mathcal S_m)$ (equivalently, $\rho_{m+1}^{-1}(\mathcal C_m)$) of
$\mathcal H_{m+1}$. We have short exact sequences
$$1\to\Ker(\rho_{m+1})\to\mathcal T_{m+1}\to\mathcal C_m\to
1\;\,\text{and}\;\,1\to\Ker(\rho_{m+1})\to\mathcal T_{m+1}^0\to\mathcal
C_m^0\to 1.$$
Thus $\dim(\mathcal T_{m+1})=\dim(\mathcal T_{m+1}^0)=r^2+\dim(\mathcal
C_m)=r^2+\dim(\mathcal C_m^0)$. 

Let $\mathcal V_{m+1}$ be the inverse image of the point $1_M[m]\in\mathcal
O_m(k)\subset\mathcal D_m(k)$ under the faithfully flat morphism $\tau_{m+1}:\mathcal
D_{m+1}\twoheadrightarrow\mathcal D_m$; it is isomorphic to an affine space over $k$
of dimension $r^2$ (cf. Subsubsection  \ref{FunctorW_m} applied to $\GL_M$). Let $\Bbb
O_{m+1,m,1_M}$ be the set of orbits of the action $\Bbb T_{m+1}$ which map to
$\mathcal O_m$ (i.e. which are contained in $\tau_{m+1}^{-1}(\mathcal O_m)$). Let
$\Bbb O_{m+1,m,1_M[m]}$ be the set of orbits of the natural action of $\mathcal
T_{m+1}$ on $\mathcal V_{m+1}$. We have the following obvious principle:

\begin{fact}\label{F2} 
There exists a natural bijection $\Bbb O_{m+1,m,1_M}\to\Bbb
O_{m+1,m,1_M[m]}$ defined by the rule: the orbit $o\in\Bbb O_{m+1,m,1_M}$ is mapped
to the element of $\Bbb O_{m+1,m,1_M[m]}$ which is the reduced scheme of
$o\cap\mathcal V_{m+1}$.
\end{fact}

\subsubsection{Geometric interpretations of sets of orbits}\label{geom}

\begin{lemma}\label{L2} The following three properties hold:

\medskip 
{\bf (a)} The set of orbits of the action $\Bbb T_m$ is in natural
bijection to the set of (representatives of) isomorphism classes of $BT_m$'s over $k$.

\smallskip
{\bf (b)} The set of orbits $\Bbb O_{m+1,m,1_M}$ (or $\Bbb O_{m+1,m,1_M[m]}$) is in
natural bijection to the set of (representatives of) isomorphism classes of
$BT_{m+1}$'s over $k$ which lift $D[p^m]$.

\smallskip
{\bf (c)} If $m\ge n_D$, then $\Bbb O_{m+1,m,1_M}$ (or $\Bbb O_{m+1,m,1_M[m]}$) has
only one orbit. 
\end{lemma}

\noindent
{\bf Proof:}
Let $\tilde B$ be a $BT_m$ over $k$. Let $\tilde D$ be a $p$-divisible group over
$k$ which lifts $\tilde B$, cf. [Il, Thm. 4.4 e)]. As $\tilde D$ has codimension $c$
and dimension $d$, its Dieudonn\'e module is isomorphic to $(M,\tilde
g\phi,\vartheta \tilde g^{-1})$ for some element $\tilde
g\in\GL_M(W(k))=\mathcal D_{\infty}(k)$. The Dieudonn\'e module of $\tilde B$
is isomorphic to $(M/p^mM,\tilde g[m]\phi_m,\vartheta_m\tilde g[m]^{-1})$. Thus (a)
and (b) are implied by Lemma \ref{L1} and the classical Dieudonn\'e theory which
provides an antiequivalence between the category of $BT_m$'s  over $k$ and the
category of Dieudonn\'e modules of $BT_m$'s over $k$ (see [Fo, pp. 153 and 160], [BBM, Chs. 3 and 4], etc.).

We prove (c). Let $B$ be a $BT_{m+1}$ over $k$ such that $B[p^m]=D[p^m]$. As $m\ge n_D$, each $p$-divisible group over $k$ that lifts $B$ is isomorphic to $D$. Thus $B$ is isomorphic to $D[p^{m+1}]$. From this and (b), we get that (c) holds.\endproof

\begin{cor}\label{cor4} The following three properties hold:

\medskip
{\bf (a)} In the situation of the last paragraph of Subsubsection \ref{limit}, each projective limit $\mathcal O_{\infty}^{\prime}:=\text{proj.}\text{lim.}_{m\to\infty}\mathcal O_m^{\prime}$ taken in the category of ringed spaces  is also a projective limit in the category of $k$-schemes and it is an integral, quasi-compact, locally closed subscheme of $\mathcal D_{\infty}$ whose $k$-valued points form one orbit of the action $\Bbb T_{\infty}(k):\mathcal H_{\infty}(k)\times\mathcal D_{\infty}(k)\to \mathcal D_{\infty}(k)$ in the category of sets. Thus $\mathcal O_{\infty}:=\text{proj.}\text{lim.}_{m\to\infty}\mathcal O_m$ is the orbit of $\Bbb T_{\infty}$ that contains $1_M\in\mathcal D_{\infty}(k)$.

\smallskip
{\bf (b)} Two points $g_1,g_2\in \mathcal D_{\infty}(k)=\GL_M(W(k))$ belong to the same orbit of the action $\Bbb T_{\infty}$ (equivalently $\Bbb T_{\infty}(k)$) if and only if the Dieudonn\'e modules
$(M,g_1\phi,\vartheta g_1^{-1})$ and  $(M,g_2\phi,\vartheta g_2^{-1})$ are
isomorphic. Therefore the set of orbits of $\Bbb T_{\infty}$ is in natural bijection
to to the set of (representatives of) isomorphism classes of $p$-divisible groups
over $k$ of codimension $c$ and dimension $d$.

\smallskip
{\bf (c)} If $K$ is an algebraically closed field that contains $k$, then $n_{D_K}=n_D$.
\end{cor}

\noindent
{\bf Proof:} The classical Dieudonn\'e theory also provides an antiequivalence between the category of $p$-divisible groups over $k$ of codimension $c$ and dimension $d$ and the category of Dieudonn\'e modules over $k$ which are isomorphic to $(M,g\phi,\vartheta g^{-1})$ for some element $g\in\GL_M(W(k))$.

Due to Lemma \ref{L2} (c), for $m\ge n_D$ we have $\mathcal O_{m+1}=\tau^{-1}_{m+1}(\mathcal O_m)$. Thus the morphism $\mathcal O_{m+1}\to\mathcal O_m$ induced by $\tau_{m+1}$ is affine (as $\tau_{m+1}$ is so). Thus $\mathcal O_{\infty}$ is an integral $k$-scheme. Moreover, $\mathcal O_{\infty}$ is the quasi-compact, locally closed subscheme of $\mathcal D_{\infty}$ which for each $m\ge n_D$ is the inverse image of $\mathcal O_{m}$ under the limit affine morphism $\tau_{\infty,m}:\mathcal D_{\infty}\to\mathcal D_m$. As $1_M\in\mathcal O_{\infty}(k)$ and $\mathcal O_{\infty}=\tau_{\infty,m}^{-1}(\mathcal O_m)$ for $m\ge n_D$, $\mathcal O_{\infty}(k)$ is the set of elements $g\in\mathcal D_{\infty}(k)=\GL_M(W(k))$ such that (cf. Lemma \ref{L2} (a)) the reduction modulo $p^m$ of $(M,g\phi,\vartheta g^{-1})$ is isomorphic to the Dieudonn\'e module of $D[p^m]$ for some $m\ge n_D$ and thus (cf. the definition of $n_D$ and the classical Dieudonn\'e theory) it is the set of elements $g\in\GL_M(W(k))$ such that there exists an isomorphism $h:(M,\phi,\vartheta)\to (M,g\phi,\vartheta g^{-1})$, i.e. we have $g=h\phi(h^{-1})$ for some element $h\in \GL_M(W(k))$. From this, Fact \ref{F1}, and the definition of $\Bbb T_{\infty}$ we get that $\mathcal O_{\infty}(k)$ is the orbit of $1_M$ under the action $\Bbb T_{\infty}(k)$. 

As $n_D$ admits upper bounds that depend only on $c$ and $d$ (cf. beginning of Section 1), there exists a $p$-divisible group $D^{\prime}$ over $k$ such that for each $m\in\Bbb N^*$ the orbit $\mathcal O_m^{\prime}$ corresponds to the isomorphism class of $D^{\prime}[p^m]$ via the bijection of Lemma \ref {L2} (a). Thus as above we argue that $\mathcal O_{\infty}^{\prime}$ is an integral,  quasi-compact, locally closed subscheme of $\mathcal D_{\infty}$ whose $k$-valued points form one orbit of the action $\Bbb T_{\infty}(k)$. Therefore part (a) holds. 

Part (b) follows from (the proof of) (a) and the classical Dieudonn\'e theory.

To prove (c), we first remark that the action $\Bbb T_{m,K}$ is the analogue of the action $\Bbb T_m$ but working with $D_K$ instead of $D$. From Lemma \ref{L2} (b) and the very definition of $n_D$ we get that: (i) we have $n_D=0$ if and only if for each $m\in\Bbb N^{\ast}$ the action $\Bbb T_m$ has only one orbit, and (ii) for $n\in\Bbb N^{\ast}$, we have $n\ge n_D$ if and only if for all $m\ge n$ the set of orbits $\Bbb O_{m+1,m,1_M}$ has only one element. From the last two sentences we get that for $n\in\Bbb N$ we have $n\ge n_D$ if and only if $n\ge n_{D_K}$. Thus $n_D=n_{D_K}$, i.e. (c) holds. 
\endproof

\medskip
 Let $\pmb{Aut}(D[p^m])_{\text{crys}}$ be the group scheme over $k$ of automorphisms of $(M/p^mM,\phi_m,\vartheta_m)$. Thus, if $R$ is a commutative $k$-algebra
and if $\sigma_R$ is the Frobenius endomorphism of the ring $W_m(R)$ of Witt vectors
of length $m$ with coefficients in $R$, then $\pmb{Aut}(D[p^m])_{\text{crys}}(R)$ is
the subgroup of $\mathcal D_m(R)=\GL_M(W_m(R))$ formed by those $W_m(R)$-linear automorphisms $\natural$ of $W_m(R)\otimes_{W_m(k)} M/p^mM$ that satisfy the identities
$(1_{W_m(R)}\otimes\phi_m)\circ\natural^{(\sigma)}=\natural\circ (1_{W_m(R)}\otimes\phi_m)$ and $\natural^{(\sigma)}\circ (1_{W_m(R)}\otimes \vartheta_m)=(1_{W_m(R)}\otimes \vartheta_m)\otimes\natural$; here $\phi_m$
and $\vartheta_m$ are viewed as $W_m(k)$-linear maps
$(M/p^mM)^{(\sigma)}\to M/p^mM$ and $M/p^mM\to (M/p^mM)^{(\sigma)}$ (respectively).
We similarly define the group scheme $\pmb{End}(D[p^m])_{\text{crys}}$ of endomorphisms of $(M/p^mM,\phi_m,\vartheta_m)$.

The following result is a version of [Va2, Thm. 2.4]. 

\begin{thm}\label{T5}
Let $m\in\Bbb N^{\ast}$. With the notations of Subsubsection
\ref{not}, the following four properties hold:

\medskip
{\bf (a)} the connected smooth affine group $\mathcal C_m^0$ is unipotent;

\smallskip
{\bf (b)} there exist two finite homomorphisms 
$$\iota_{m}:\mathcal
S_m\to\pmb{Aut}(D[p^m])_{\text{crys}}\;\;\;\text{and}\;\;\;\zeta_{m}:\pmb{Aut}(D[p^m])\to
\pmb{Aut}(D[p^m])_{\text{crys}}$$ which at the level of $k$-valued points
induce isomorphisms $\iota_m(k):\mathcal
S_m(k)\to\pmb{Aut}(D[p^m])_{\text{crys}}(k)$ and
$\zeta_{m}(k):\pmb{Aut}(D[p^m])(k)\to
\pmb{Aut}(D[p^m])_{\text{crys}}(k)$;

\smallskip
{\bf (c)} we have $\dim(\mathcal S_m)=\dim(\mathcal C_m)=\dim(\mathcal
C^0_m)=\gamma_D(m)$;

\smallskip
{\bf (d)} the connected smooth affine group $\mathcal T_{m+1}^0$ is unipotent and has dimension $r^2+\gamma_D(m)$.
\end{thm}

\noindent
{\bf Proof:}
Parts (a) to (c) are proved in [Va2, Thm. 2.4] using reduced groups. We will only recall here the
definitions of $\zeta_m$ and $\iota_m$. We will view $\mathcal D_m$
as a connected smooth affine group over $k$. 

The crystalline Dieudonn\'e theory provides us with a homomorphism 
$$\zeta_m:\pmb{Aut}(D[p^m])\to \pmb{Aut}(D[p^m])_{\text{crys}}$$
that maps a point $x\in \pmb{Aut}(D[p^m])(R)$ to the inverse of the element of
$\pmb{Aut}(D[p^m])_{\text{crys}}(R)$ which is the evaluation of $\Bbb D(x)$ at the
thickening $\Spec R \hookrightarrow\Spec W_m(R)$ defined by the natural divided
power structure on the kernel of the epimorphism $W_m(R)\twoheadrightarrow R$. 

The homomorphism $\tilde{\mathcal P}:\tilde{\mathcal H}\to \pmb{GL}_M$ provides us with a homomorphism
$$\iota_m:\mathcal S_m\to \pmb{Aut}(D[p^m])_{\text{crys}}$$ 
that maps a point $\tilde h_R[m]\in \mathcal S_m(R)\subset
\mathcal H_m(R)=\tilde{\mathcal H}(W_m(R))$ to the point $h_R[m]:=\tilde{\mathcal P}(W_m(R))(\tilde h_R[m])\in\pmb{Aut}(D[p^m])_{\text{crys}}(R)\subset \mathcal D_m(R)=\GL_M(W_m(R))$.

We prove (d). The group $\Ker(\rho_{m+1})$ is a unipotent, commutative group of
dimension $r^2$, cf. Subsubsection  \ref{FunctorW_m}. As the connected group $\mathcal C_m^0$
is  unipotent (cf. (a)) of dimension $\gamma_D(m)$ (cf. (c)) and as we have a short
exact sequence $1\to \Ker(\rho_{m+1})\to\mathcal T_{m+1}^0\to\mathcal C_m^0\to 1$
(see Subsubsection \ref{not}), from [DGr, Vol. II, Exp. XVII, Prop. 2.2 iv)] we get that
the group $\mathcal T_{m+1}^0$ is a unipotent group of dimension $r^2+\gamma_D(m)$. Thus (d) holds.
\endproof

\begin{cor}\label{cor5}
 The smooth groups $\mathcal C_m$, $\pmb{Aut}(D[p^m])_{\text{crys},\re}$, and
$\pmb{Aut}(D[p^m])_{\re}$ are radicially isogenous.
\end{cor}

\noindent
{\bf Proof:}
The finite homomorphisms $\iota_m$ and $\zeta_m$ of Theorem \ref{T5} (b) are bijective. Thus the corollary follows from [Se, Subsect. 1.2, Props. 1 and 2]
(loc. cit. is stated for commutative group schemes over $k$ but its proofs do not 
require the commutativity assumption).
\endproof

\subsubsection{Two formulas}\label{2flas}
As $\gamma_D(m)=\dim(\mathcal S_m)$ (cf. Theorem \ref{T5} (d)), we have 
$$\gamma_D(m)=\dim(\mathcal H_m)-\dim(\mathcal O_m). \leqno (1) $$
As $\dim(\mathcal H_{m+1})-\dim(\mathcal H_m)=r^2$ we get that
$$\gamma_D(m+1)-\gamma_D(m)=r^2-\dim(\mathcal O_{m+1})+\dim(\mathcal O_m).\leqno (2)$$

\subsubsection{Digression: a variant of the action $\Bbb T_m$}\label{digr}
For the sake of completeness, we describe here a variant of the action $\Bbb T_m$
that appeals to Subsubsection \ref{tildeH}. 

We consider the functor $\text{Aff}_k\to \text{Sets}$ which takes $\Spec R$ to
$W_m(R)$-linear isomorphisms from $W_m(R)\otimes_{W(k)}M^{(\sigma)}$ to
$W_m(R)\otimes_{W(k)} M$. This functor is represented by a connected smooth affine $k$-scheme of finite type $\mathcal D^\prime_m$. If $\mathcal
H^\prime_m:=\Bbb W_m(\pmb{Aut}(\phi,\vartheta))$, then we have a natural action
$\Bbb T_m^\prime:\mathcal H_m^\prime\times_k \mathcal D_m^\prime\to \mathcal
D_m^\prime$ defined by transport of structure. 

One can check that $\Bbb T_m^\prime$ is isomorphic to $\Bbb
T_m^{(\sigma)}$ and has properties similar to the ones enjoyed by $\Bbb T_m$ (see below; the isomorphism $\mathcal H_m^\prime=\Bbb W_m(\pmb{Aut}(\phi,\vartheta))\to \Bbb W_m(\tilde{\mathcal H}^{(\sigma)})=\mathcal H_m^{(\sigma)}$ is $\Bbb W_m$ of the isomorphism $\pmb{Aut}(\phi,\vartheta)\to \tilde{\mathcal H}^{(\sigma)}$ of Subsubsection \ref{tildeH}). In particular, one can use the action $\Bbb T_m^\prime$ to show
that the three smooth groups of the Corollary \ref{cor5} are in fact isomorphic. More precisely, $\iota_m$ can be identified with the Frobenius endomorphism of $\mathcal{S}_m$, $\zeta_{m,\re}$ can be identified with the $m^{\text{th}}$-power of the Frobenius endomorphism of $\pmb{Aut}(D[p^m])_{\re}$, and $\Ker(\zeta_m)$ contains (and it is quite likely to be equal to) the kernel of the $m^{\text{th}}$-power of the Frobenius endomorphism of $\pmb{Aut}(D[p^m])$.

Let $S=\Spec A$ be a smooth affine $k$-scheme. Let the crystalline site $\text{CRIS}(S/W_m(k))$ be as in
[BBM, Subsect. 1.1.1] over the natural divided power thickening
$\Spec k\hookrightarrow \Spec W_m(k)$. One identifies the
sections of $\mathcal O_{\text{CRIS}(S/W_m(k))}$ with
$W_m(A^{(\sigma^m)})$, cf. [IR, III, (1.5.2), p. 142]. 

For $i\in\{1,2\}$ let $B_i$ be a $BT_m$ over $k$
with Dieudonn\'e $F$-crystal $\Bbb D(B_i)=(N_i,\phi_{m,i},\vartheta_{m,i})$. Due to the previous paragraph, morphisms $\Bbb D(B_{2,S})\to \Bbb D(B_{1,S})$ of
Dieudonn\'e  $F$-crystals of $\text{CRIS}(S/W_m(k))$ correspond to $W_m(A^{(\sigma^m)})$-linear maps
$W_m(A^{(\sigma^m)})\otimes_{W_m(k)} N_2\to W_m(A^{(\sigma^m)})\otimes_{W_m(k)} N_1$ compatible with the $\phi_{m,i}$'s and $\vartheta_{m,i}$'s. In view of
[BM, Thm. 4.1.1], such $W_m(A^{(\sigma^m)})$-linear maps describe the abelian
group $\Hom(B_{1,S},B_{2,S})$. One can check via a direct inspection that the orbits
of the action $\Bbb T_m^\prime$ are in bijection to the set of (representatives) of
isomorphism classes of $BT_m$'s over $k$, to be compared with [Il, proof b) of Prop. 1.7].

\subsection{Proof of Theorem \ref{UP}}\label{PTUP} 
The set of (representatives of) isomorphism classes of $BT_{m+1}$'s over $k$ which
lift $D[p^m]$ is in bijection to the set $\Bbb O_{m+1,m,1_M[m]}$ of orbits of the
action of $\mathcal T_{m+1}$ on $\mathcal V_{m+1}$, cf. Lemma \ref{L2} (b).
The group $\mathcal T_{m+1}^0$ is unipotent of dimension $r^2+\gamma_D(m)$, cf.
Theorem \ref{T5} (d). Moreover $\mathcal V_{m+1}$ is isomorphic to the affine
space of dimension $r^2$ over $k$. We have a short exact sequence
$1\to\Ker(\rho_{m+1})\to\mathcal T_{m+1}\to\mathcal C_m\to 1$, where
$\Ker(\rho_{m+1})$ is isomorphic to $\Bbb G_a^{r^2}$ (cf. Subsubsection  \ref{FunctorW_m}) and
where $\mathcal C_m$ is radicially isogenous to $\pmb{Aut}(D[p^m])_{\re}$ (cf.
Corollary \ref{cor5}). Theorem \ref{UP} follows from the last four
sentences.\endproof 

\subsection{Properties of orbits}\label{orbits}

\begin{prop}\label{P1}
Let $U$ be a smooth affine group over $k$ whose identity component $U^0$ is
unipotent. Let $V$ be a reduced affine variety over $k$ equipped with an action
$\Lambda:U\times_k V\to V$. Then all orbits of $\Lambda$ are closed. Thus, if
$\Lambda$ has an open, Zariski dense orbit $o$, then $\Lambda$ has a unique orbit (i.e.
$o=V$).
\end{prop}

\noindent
{\bf Proof:}
All orbits of the action of $U^0$ on $V$ are closed, cf. [Sp, Prop. 2.4.14] or [DGa, Ch. IV, Sect. 2, Subsect. 2, Cor. 2.7]. Each orbit of $\Lambda$ is set-theoretically a finite union of orbits of the action of $U^0$ on $V$. Thus, the orbits of $\Lambda$ are closed. The last part of the proposition follows from the fact that $V$ is the only closed, Zariski dense subscheme of itself.
\endproof

\medskip
We will need the following three applications of Proposition \ref{P1} in the proofs of Theorems \ref{Gamma} and \ref{INF}.

\begin{cor}\label{cor6}
The following three properties hold:

\medskip
{\bf (a)} The open embeddings $\pmb{Aut}(D[p^m])\hookrightarrow \pmb{End}(D[p^m])$ and $\pmb{Aut}(D[p^m])_{\text{crys}}\hookrightarrow \pmb{End}(D[p^m])_{\text{crys}}$ are also closed.

\smallskip
{\bf (b)} The identity component $\pmb{Aut}(D[p^m])^0$ (resp. $\pmb{Aut}(D[p^m])^0_{\text{crys}}$) is the translation of the identity component $\pmb{End}(D[p^m])^0$ (resp. $\pmb{End}(D[p^m])^0_{\text{crys}}$) via the $k$-valued point $1_{D[p^m]}$ (resp. $1_{M/p^mM}$). 

\smallskip
{\bf (c)} Each $k$-valued point of $\pmb{End}(D[p^m])^0_{\text{crys}}$ is a nilpotent endomorphism of $M/p^mM$.
\end{cor}

\noindent
{\bf Proof:} We will only prove (a) and (b) for $\pmb{Aut}(D[p^m])$ and $\pmb{End}(D[p^m])$ as the arguments for their crystalline versions are the same. To prove (a), it suffices to show that $U:=\pmb{Aut}(D[p^m])_{\re}$ is an open and closed subvariety of $V:=\pmb{End}(D[p^m])_{\re}$. As the identity component of $U$ is a unipotent group (cf. Theorem \ref{T5} (a) and (b)) and as $U$ is an orbit of the canonical left composite action of $U$ on the affine variety $V$, from Proposition \ref{P1} we get that $U$ is an open and closed subvariety of $V$. Part (b) follows from (a).

To prove (c), we can assume that $m=1$. As each $k$-valued point of the unipotent subgroup $\pmb{Aut}(D[p])^0_{\text{crys},\re}$ of $\GL_{M/pM}$ is unipotent, the case $m=1$ follows from (b).

For the sake of completeness we include an elementary second proof of the corollary that does not depend on Proposition \ref{P1}. It suffices to prove that if $\mathcal U$ is a connected component of $\pmb{End}(D[p])_{\text{crys}}$, then the characteristic polynomial $\chi_e(T)=\det(T1_M[1]-e)$ of an element $e\in\mathcal U(k)$  belongs to $\Bbb F_p[T]$ and thus due to continuity reasons it depends only on $\mathcal U$  (so if $\mathcal U=\pmb{End}(D[p])_{\text{crys}}^0$, then $\chi_e(T)=\chi_0(T)=T^r$). As $e$ commutes with $\phi_1$ and $\vartheta_1$ and as we have $\Image(\phi_1)=\Ker(\vartheta_1)$ and $\Ker(\phi_1)=\Image(\vartheta_1)$, the short exact sequences of $k$-vector spaces $0\to \Image(\phi_1)\to M/pM\to (M/pM)/\Image(\phi_1)\to 0$ and $0\to (M/pM)/\Image(\phi_1)\to (M/pM)^{(\sigma)}\to  \Image(\phi_1)\to 0$ are invariant under $e$ and $e^{(\sigma)}$. Thus $\chi_e(T)=\chi_{e^{(\sigma)}}(T)$ and therefore $\chi_e(T)\in\Bbb F_p[T]$.\endproof

\begin{lemma}\label{L3}
The following three statements are equivalent:

\medskip
{\bf (i)} we have $\dim(\mathcal O_{m+1})=\dim(\mathcal
O_m)+\dim(\Ker(\rho_{m+1}))=\dim(\mathcal O_m)+r^2$;

\smallskip
{\bf (ii)} we have $\mathcal O_{m+1}=\tau_{m+1}^{-1}(\mathcal O_m)$ (i.e. we have
$\mathcal V_{m+1}\subset \mathcal O_{m+1}$);

\smallskip
{\bf (iii)} the set of orbits $\Bbb O_{m+1,m,1_M}$ (equivalently, $\Bbb
O_{m+1,m,1_M[m]}$) is finite.
\end{lemma}

\noindent
{\bf Proof:}
The $k$-scheme $\mathcal V_{m+1}$ is an affine space of dimension $r^2$. The fibres
of the faithfully flat morphism $\mathcal O_{m+1}\to\mathcal O_m$ induced by $\tau_{m+1}$ have equal
dimensions. Thus $\dim(\mathcal O_{m+1})=\dim(\mathcal O_m)+\dim(\mathcal
V_{m+1}\cap\mathcal O_{m+1})$. If (ii) holds, then $\dim(\mathcal
V_{m+1}\cap\mathcal O_{m+1})=\dim(\mathcal V_{m+1})=r^2$ and therefore
$\dim(\mathcal O_{m+1})=\dim(\mathcal O_m)+r^2$. Thus (ii) implies (i). Obviously
(ii) implies (iii). 

We check that (i) implies (ii). The orbit of $1_M[m+1]$ under $\mathcal T_{m+1}$ is
$(\mathcal V_{m+1}\cap\mathcal O_{m+1})_{\re}$. If (i) holds, then
$\dim(\mathcal V_{m+1}\cap\mathcal O_{m+1})=\dim(\mathcal O_{m+1})-\dim(\mathcal
O_m)=r^2=\dim(\mathcal V_{m+1})$ and therefore $\mathcal V_{m+1}\cap\mathcal
O_{m+1}$ is an open subscheme of $\mathcal V_{m+1}$. Thus the action of $\mathcal
T_{m+1}$ on $\mathcal V_{m+1}$ has an open, Zariski dense  orbit. As $\mathcal T_{m+1}^0$ is
a unipotent group (cf. Theorem \ref{T5} (d)), from the second part of Proposition
\ref{P1} we get that $\mathcal V_{m+1}\cap\mathcal O_{m+1}=\mathcal V_{m+1}$, i.e.
$\mathcal V_{m+1}\subset \mathcal O_{m+1}$. Thus (i) implies (ii). 

We check that (iii) implies (ii). As (iii) holds, the action of $\mathcal T_{m+1}$
on $\mathcal V_{m+1}$ has a finite number of orbits and thus it has an open, Zariski dense
orbit. As in the previous paragraph we argue that this orbit is $\mathcal V_{m+1}$, i.e. $\mathcal V_{m+1}\subset \mathcal O_{m+1}$.
Therefore (iii) implies (ii). Thus the three statements are equivalent.
\endproof

\begin{lemma}\label{L4}
If $m\in\{1,\ldots,n_D-1\}$ (thus $n_D\ge 2$), then $\gamma_D(m+1)\neq \gamma_D(m)$.
\end{lemma}

\noindent
{\bf Proof:} It suffices to show
that the assumptions that $m\in\{1,\ldots,n_D-1\}$ and $\gamma_D(m+1)=\gamma_D(m)$ lead to a contradiction. 
As $\gamma_D(m+1)-\gamma_D(m)=0$, from Theorem \ref{Gamma} (a) (proved in Subsection \ref{SSect2.2}) we get that the sequence $(\gamma_D(i))_{i\ge m}$ is constant, i.e. for each
integer $i\ge m$ we have $\gamma_D(i)=\gamma_D(m)$. In particular, as
$m\in\{1,\ldots,n_D-1\}$ we get that $\gamma_D(n_D)=\gamma_D(n_D-1)$.

Thus we have $\dim(\mathcal O_{n_D})=r^2+\dim(\mathcal O_{n_D-1})$, cf. Formula (2).
Therefore the set of orbits $\Bbb O_{n_D,n_D-1,1_M}$ (equivalently, $\Bbb
O_{n_D,n_D-1,1_M[n_D-1]}$) has one element (cf. the equivalence of statements (i) and (ii) of Lemma \ref{L3}). Thus each $BT_{n_D}$ over $k$ that lifts $D[p^{n_D-1}]$
is isomorphic to $D[p^{n_D}]$, cf. Lemma \ref{L2} (b). 

Let $D^\prime$ be a $p$-divisible group over $k$ such that $D^\prime[p^{n_D-1}]$ is isomorphic
to $D[p^{n_D-1}]$. We know that $D^\prime[p^{n_D}]$ is isomorphic to $D[p^{n_D}]$, cf.
previous paragraph. Thus $D^\prime$ is isomorphic to $D$, cf. the definition of
$n_D$. Therefore we have $n_D\le n_D-1$ (again cf. the definition of $n_D$), a
contradiction.\endproof

 \subsection{Proof of Theorem \ref{Gamma}}\label{PGamma} 
We know that Theorem \ref{Gamma} (a) holds, cf. Subsection \ref{SSect2.2}.
From Theorem \ref{Gamma} (a) and Lemma \ref{L4} we get that the strict inequalities of Theorem \ref{Gamma} (b) hold. Thus to end the proof of Theorem \ref{Gamma} (b) it suffices to show that $\gamma_D(1)\ge a_D^2$. Taking the Cartier dual of the exact complex $0\to \pmb{\alpha}_p^{a_{D^{\text{t}}}}\to D[p]^{\text{t}}$, we get an exact complex $D[p]\to \pmb{\alpha}_p^{a_{D^{\text{t}}}}\to 0$. Using it and the exact complex $0\to \pmb{\alpha}_p^{a_D}\to D[p]$, we get that the group scheme of homomorphisms $\pmb{Hom}(\pmb{\alpha}_p^{a_{D^{\text{t}}}},\pmb{\alpha}_p^{a_D})=\pmb{End}(\pmb{\alpha}_p)^{a_Da_{D^{\text{t}}}}$ is a subgroup scheme of $\pmb{End}(D[p])$. At the level of dimensions we get that $\gamma_D(1)\ge a_Da_{D^{\text{t}}}$. Thus to prove that $\gamma_D(1)\ge a_D^2$ it suffices to check that $a_D=a_{D^{\text{t}}}$. We will use the classical Dieudonn\'e theory. We recall that the Dieudonn\'e module of $D[p]$ is $(M/pM,\phi_1,\vartheta_1)$. Then $a_D=\dim_k((M/pM)/(\Image(\phi_1)+\Image(\vartheta_1)))$ is equal to $r-\dim_k(\Image(\phi_1))-\dim_k(\Image(\vartheta_1))+\dim_k(\Image(\phi_1)\cap\Image(\vartheta_1))$ and thus also to $r-c-d+\dim_k(\Ker(\vartheta_1)\cap\Ker(\phi_1))=\dim_k(\Ker(\vartheta_1)\cap\Ker(\phi_1))=a_{D^{\text{t}}}$.

This ends the proof of Theorem \ref{Gamma} (b) and thus also of Theorem
\ref{Gamma} (cf. Subsection \ref{SSect2.2} and paragraph after Theorem \ref{Gamma}). For the sake of completeness, here is a self-contained proof of Theorem \ref{Gamma} (c) that does not rely on [Va2]. 

To prove Theorem \ref{Gamma} (c) we can assume that $n_D\ge 1$. For $i\in\Bbb N^{\ast}$ we have $\gamma_D(i)\le\dim_k(\Lie({\bf Aut}(D[p^i])))$. As the $k$-vector space $\Lie({\bf Aut}(D[p^i]))$ is isomorphic to $t_{D[p]}\otimes_k t_{D[p]^{\text{t}}}$ (cf. [Il, Prop. 3.1 b) and Lemma 4.1]; here $t_{D[p]}$ is the zero cohomology group of the Lie complex of $D[p]$ as used in [Il, Sect. 2.1]),  it has dimension $cd$. Thus $\gamma_D(i)\le cd$. If moreover $i\ge n_D$, then from  Lemma \ref{L2} (c) we get that we have $\mathcal O_{i+1}=\tau_{i+1}^{-1}(\mathcal O_i)$ and thus from the equivalence of statements (i) and (ii) of Lemma \ref{L3} we get that $\dim(\mathcal O_{i+1})=\dim(\mathcal O_i)+r^2$ and therefore that (cf. Formula (2)) $\gamma_D(i+1)=\gamma_D(i)$. We conclude that Theorem \ref{Gamma} (c) holds.\endproof

\subsection{Proof of Theorem \ref{CP}}\label{PCP} 
We prove that Theorem \ref{CP} holds. The
integral domain $\mathcal R$ is normal and its field of fractions $k_1$ is
algebraically closed. Thus $\mathcal R$ is a perfect ring. Let $\sigma_{\mathcal R}$
be the Frobenius automorphism of the ring $W(\mathcal R)$ of ($p$-typical) Witt vectors with
coefficients in $\mathcal R$. Let $(M_{\mathcal R},\phi_{\mathcal R},\vartheta_\mathcal
R)$ be the Dieudonn\'e module of $\mathfrak{D}$ over $\mathcal R$ (i.e. the
projective limit indexed by $s\in\Bbb N^{\ast}$ of the evaluation of the Dieudonn\'e
$F$-crystal $\Bbb D(\mathfrak{D})$ of $\mathfrak{D}$ at the thickening
$\Spec \mathcal R\hookrightarrow\Spec W_s(\mathcal R)$ defined by the natural
divided power structure of the ideal $(p)$ of $W_s(\mathcal R)$). Thus $M_\mathcal
R$ is a free $W(\mathcal R)$-module of rank $r$, $\phi_{\mathcal R}:M_{\mathcal R}\to
M_{\mathcal R}$ is a $\sigma_{\mathcal R}$-linear endomorphism, $\vartheta_{\mathcal
R}:M_{\mathcal R}\to M_{\mathcal R}$ is a $\sigma_{\mathcal R}^{-1}$-linear endomorphism, and we have $\phi_{\mathcal R}\circ \vartheta_{\mathcal R}=\vartheta_{\mathcal R}\circ\phi_{\mathcal R}=p 1_{M_{\mathcal R}}$. Let $\bar F^1_{\mathcal R}$ be the
kernel of the reduction modulo $p$ of $\phi_{\mathcal R}$. As $\mathfrak{D}$ has
codimension $c$ and dimension $d$ and as $\mathcal R$ is local, the pair
$(M_{\mathcal R},\bar F^1_{\mathcal R}\subset M_{\mathcal R}/pM_{\mathcal R})$ is
isomorphic to $W(\mathcal R)\otimes_{W(k)} (M,\bar F^1\subset M/pM)$. Thus there
exists an element $g_{\mathcal R}\in\GL_M(W(\mathcal R))$ such that the triple
$(M_{\mathcal R},\phi_{\mathcal R},\vartheta_{\mathcal R})$ is isomorphic to the triple
$(W(\mathcal R)\otimes_{W(k)} M,g_{\mathcal R}(\sigma_\mathcal
R\otimes\phi),(\sigma_{\mathcal R}^{-1}\otimes\vartheta)g_{\mathcal R}^{-1})$. 
We view $g_{\mathcal R}$ as a point $g_{\mathcal R}\in\mathcal D_{\infty}(\mathcal R)$.
For $m\in\Bbb N^{\ast}$, let $g_{\mathcal R}[m]\in\mathcal D_m(\mathcal R)$ be the
point defined by the reduction modulo $p^m$ of $g_{\mathcal R}$. Let
$g_{k_1}[m]\in\mathcal D_m(k_1)$ and $g_k[m]\in\mathcal D_m(k)$ be the points
defined naturally by $g_{\mathcal R}[m]$. 

Let $g^\prime\in\GL_M(W(k))$ be such that the Dieudonn\'e module of $D^\prime$
is isomorphic to $(M,g^\prime\phi,\vartheta (g^{\prime})^{-1})$. Let $\mathcal O_m^\prime$ be the orbit of $g^\prime[m]\in\mathcal D_m(k)$
under the action $\Bbb T_m$. As $\mathfrak{D}_{k_1}$ is isomorphic to $D^\prime_{k_1}$, for all $m\in\Bbb N^{\ast}$ we have
$g_{k_1}[m]\in\mathcal O_m^\prime(k_1)\subset\mathcal D_m(k_1)$ (cf. Lemma
\ref{L1} applied to $\Bbb T_{m,k_1}$). As $D=\mathfrak{D}_k$, we have $g_k[m]\in\mathcal O_m(k)$.
Thus, as $g_{k_1}[m]$ specializes to $g_k[m]$ (via $g_{\mathcal R}[m]$), the schematic closure $\bar{\mathcal O}_m^\prime$ of $\mathcal O_m^\prime$
in $\mathcal D_m$ contains points of $\mathcal O_m$. Thus $\mathcal O_m\subset
\bar{\mathcal O}_m^\prime$ and therefore we have
$$\dim(\mathcal O_m)\le\dim(\mathcal O_m^\prime)=\dim(\bar{\mathcal
O}_m^\prime).\leqno (3)$$

Similarly to Formula (1) we have $\gamma_{D^\prime}(m)=\dim(\mathcal H_m)-\dim(\mathcal O^\prime_m)$. Thus
for all $m\in\Bbb N^{\ast}$ we have $\gamma_{D^\prime}(m)\le\gamma_D(m)$, cf. either (3) or fact (a) of Subsection \ref{SSect1.1}.

We assume that there exists a number $l\in\Bbb N^{\ast}$ such that
$\gamma_{D^\prime}(l)=\gamma_D(l)$. Thus $\dim(\mathcal O_l)=\dim(\mathcal O_l^\prime)=\dim(\bar{\mathcal O}_l^\prime)$, cf. Formula (1). From these equalities and the fact that $\bar{\mathcal O}_l^\prime$ contains
points of $\mathcal O_l$, we get that the two orbits $\mathcal O_l^\prime$ and $\mathcal O_l$ coincide. Thus $D^\prime[p^l]$ is isomorphic to $D[p^l]$, cf. Lemma
\ref{L2} (a). This ends the proof of
Theorem \ref{CP} (a).

We prove Theorem \ref{CP} (b). Let $m\ge\max\{n_D,n_{D^\prime}\}$. We have
$s_{D^\prime}=\gamma_{D^\prime}(m)\le\gamma_D(m)=s_D$. Next we assume that $s_{D^\prime}=s_D$. Thus
$\gamma_{D^\prime}(m)=\gamma_D(m)$. Therefore $D^\prime[p^m]$ is isomorphic to $D[p^m]$,
cf. Theorem \ref{CP} (a). As $m\ge n_D$, we get that $D^\prime$ is isomorphic to
$D$. This ends the proof of Theorem \ref{CP} (b) and thus also of Theorem
\ref{CP}.\endproof

\subsection{Proof of Theorem \ref{INF}}\label{PINF}

Let $B$, $D$, and $1\le m<n_D$ be as in Theorem \ref{INF}. To prove Theorem \ref{INF} it
suffices to show that there exists an infinite set $\mathcal I$ of $BT_{m+1}$'s over
$k$ that has the following two properties:

\medskip
{\bf (i)} if $B_1$, $B_2\in \mathcal I$ are distinct elements, then $B_1$ is not
isomorphic to $B_2$; 

\smallskip
{\bf (ii)} if $B_1\in\mathcal I$, then $B_1[p^m]$ is isomorphic to $B=D[p^m]$. 

\medskip\noindent
Based on the first sentence of Subsection \ref{PTUP}, the existence of $\mathcal
I$ is equivalent to the fact that the set of orbits $\Bbb O_{m+1,m,1_M[m]}$
introduced in Subsubsection \ref{not} is infinite. We will show that the assumption that
the set of orbits $\Bbb O_{m+1,m,1_M[m]}$ is finite, leads to a contradiction. Based
on Lemma \ref{L3}, we get that $\dim(\mathcal O_{m+1})=\dim(\mathcal O_m)+r^2$. 
Thus $\gamma_D(m+1)-\gamma_D(m)=0$ (cf. Formula (2)) and this contradicts Lemma \ref{L4}.\endproof 

\section{Some applications of Theorem \ref{Gamma}}\label{Sect5}
We begin this section by first proving the corollaries of Section 1.

\subsection{Proof of Corollary \ref{cor1}}\label{Pcor1} 
To prove Corollary \ref{cor1} we can
assume that $cd>0$. If $D$ is ordinary, then we have $\gamma_D(i)=0$ for all $i\in\Bbb N$ and $n_D=1\le cd$. Thus we can also assume that $D$ is not ordinary (equivalently, that $a_D\ge 1$). From Theorem \ref{Gamma} (b) and (c) we get that $n_D\le cd+1-\gamma_D(1)\le cd+1-a_D^2\le cd$. Thus Corollary \ref{cor1} holds.\endproof

\subsection{Proof of Corollary \ref{cor2}}\label{Pcor2}
 To prove Corollary \ref{cor2} (a) we can
assume $cd>0$. Each connected component of the image of $r_{l,i}$ has dimension
$\gamma_D(l)-\gamma_D(l-i)$, cf. property \ref{SSect2.1} (iii). If an endomorphism of $D[p^l]$
restricts to an automorphism of $D[p^i]$, then it is an automorphism. Thus the image
of the restriction homomorphism $a_{l,i}:\pmb{Aut}(D[p^l])\to \pmb{Aut}(D[p^i])$ is
$\Image(r_{l,i})\cap \pmb{Aut}(D[p^i])$ and therefore it is a nonempty open
subscheme of $\Image(r_{l,i})$. Thus the images of $a_{l,i}$ and $r_{l,i}$ have
equal dimension $\gamma_D(l)-\gamma_D(l-i)$. We have $\gamma_D(l)-\gamma_D(l-i)=0$
if and only if $l-i\ge n_D$, cf. Theorem \ref{Gamma} (b) and (c). Therefore the
image of $a_{l,i}$ (or $r_{l,i}$) is finite if and only if $l-i\ge n_D$. Thus Corollary \ref{cor2} (a) holds.

From Corollary \ref{cor2} (a) and Theorem \ref{Gamma} (b) and (c) we get that if $D$
is not ordinary (i.e. if $n_D>0$ and $\gamma_D(1)>0$) and if $n$ is as in Corollary
\ref{cor2} (b), then we have $n=n_D$. Thus Corollary \ref{cor2} (b) holds.

To prove Corollary \ref{cor2} (c), we can assume that $D$ is not ordinary. Thus Corollary \ref{cor2} (c) follows from Corollary \ref{cor2} (b) and the fact that the image of the restriction homomorphism $r_{l,1}$ is finite if and only if the image of the analogous restriction homomorphism obtained working with $D^s$ instead of $D$ is finite.\endproof

\subsection{Proof of Corollary \ref{cor3}}\label{Pcor3}
Corollary \ref{cor3} (a) follows from Corollary \ref{cor2} (a). Corollary \ref{cor3} (b) follows from Theorem \ref{T5} (a) and (b). To check that Corollary \ref{cor3}
(c) to (e) hold, we can assume that $D$ is not ordinary and that $i\ge n_D+1\geq
2$. Let $a_i\in\pmb{Aut}(D[p^i])^0(k)$. Let $b_i:=\zeta_i(k)(a_i)\in
\mathcal D_i(k)=\GL(W_i(k))$ be the crystalline realization of $a_i$, cf. the
notations of Subsection \ref{act}. Due to Corollary \ref{cor3} (a), we can write
$b_i=1_M[i]+p^{s(a_i)}e_i$, where $s(a_i)\in\{i-n_D,\ldots,i\}$ and
$e_i\in\End(M/p^iM)\setminus p\End(M/p^iM)$; moreover there exist elements $a_i$ with $s(a_i)=i-n_D$. Corollary \ref{cor3} (c) follows from this once we check that the order of the element $b_i$ (equivalently $a_i$) is $p^{i-s(a_i)}$. It suffices to consider the case when $n_D\ge 2$ and $(p,i,s(a_i))=(2,n_D+1,1)$ and to show that $b_i^2=1_M[i]+4(e_i+e_i^2)$ is not congruent to $1_M[i]$ modulo $8$. In such a case, the reduction $\bar e_i$ modulo $2$ of $e_i$ is a non-zero nilpotent element of $\End(M/2M)$ (the elements $\bar e_i$ together with $0$ are the $k$-valued points of a connected smooth subgroup of $\pmb{End}(D[2])^0_{\text{crys}}$ and thus are nilpotent, cf. Corollary \ref{cor6} (c) applied with $m=1$); thus $\bar e_i\neq \bar e_i^2$ (i.e. $b_i^2$ is not congruent to $1_M[i]$ modulo $8$). 

If $a_i'\in\pmb{Aut}(D[p^i])^0(k)$ is another element whose crystalline realization is 
$b_i'=1_M[i]+p^{s(a_i')}e_i'\in\GL(W_i(k))$ with $e_i'\in\End(M/p^iM)\setminus p\End(M/p^iM)$, then $b_i$ and $b_i'$ commute
modulo $p^{\min\{s(a_i)+s(a_i'),i\}}$. If $i\ge 2n_D$, then $s(a_i)+s(a_i')\geq
2(i-n_D)\geq i$ and thus $b_i$ and $b_i'$ commute. Thus Corollary
\ref{cor3} (d) holds. If $l>i\ge n_D+1$ and $a_l\in\pmb{Aut}(D[p^l])^0(k)$ has crystalline
realization $b_l=1_M[l]+p^{s(a_l)}e_l$, then there exists a unique element
${}_ia_l\in\pmb{Aut}(D[p^i])^0(k)$ whose crystalline realization is
${}_ib_l=1_M[i]+p^{s(a_l)-l+i}e_l[i]$, where $e_l[i]$ is the reduction modulo $p^i$
of $e_l$. The rule $a_l\mapsto {}_ia_l$ defines the canonical bijection
$\pmb{Aut}(D[p^l])^0(k)\to \pmb{Aut}(D[p^i])^0(k)$ which has the desired properties. Thus Corollary \ref{cor3} (e) holds.\endproof

\subsection{A refinement of Corollary \ref{cor1}}\label{Rcor1}
 We recall that
$s_D=\gamma_D(n_D)$. If $D$ is ordinary, then $s_D=0=a_D$ and $n_D\le 1$. Thus in this
subsection we will assume that $D$ is not ordinary. From Theorem \ref{Gamma} (b) we
get that $n_D\le s_D+1-\gamma_D(1)\le s_D+1-a_D^2$.
As $D$ is not ordinary, we have $a_D\ge 1$ and
thus we
have 
$$n_D\le s_D+1-a_D^2\le s_D.\leqno (4)$$
\indent
Dieudonn\'e's
classification of $F$-isocrystals over $k$ implies that we have a direct sum
decomposition $(M[{1\over p}],\phi)=\oplus_{s=1}^v (W_s,\phi)$
into simple $F$-isocrystals over $k$ (here $v$ is a positive integer). More
precisely, for $s\in\{1,\ldots,v\}$ there exist $c_s,d_s\in\Bbb N$ such that
$r_s:=c_s+d_s>0$, $g.c.d.\{c_s,d_s\}=1$, $\dim_{B(k)}(W_s)=r_s$, and moreover there
exists a $B(k)$-basis for $W_s$ formed by elements fixed by $p^{-d_s}\phi^{r_s}$;
the unique Newton polygon slope of $(W_s,\phi)$ is $\alpha_s:={{d_s}\over {r_s}}\in
\Bbb Q\cap [0,1]$. We have 
$$s_D=cd-{1\over 2}\sum_{s=1}^v\sum_{t=1}^v r_sr_t|\alpha_s-\alpha_t|=cd-{1\over
2}\sum_{s=1}^v\sum_{t=1}^v |c_sd_t-c_td_s|,\leqno (5) $$
cf. [Va2, Thm. 1.2 (c), (d), (e), and (f)]. From (4) and (5) we get 
$$n_D\le 1+cd-a_D^2-{1\over 2}\sum_{s=1}^v\sum_{t=1}^v |c_sd_t-c_td_s|\le cd-{1\over 2}\sum_{s=1}^v\sum_{t=1}^v |c_sd_t-c_td_s|.\leqno (6)$$

\subsection{Remark}\label{Remark1} 
In the proof of [Va2, Thm. 1.2 (e)] (see [Va2, Subsect. 3.5])
one needs to take $m$ slightly bigger (it would suffice to take $m\ge
\max\{2\kappa+n_{\tilde D},n_D\}$ instead of $m\ge \max\{\kappa+n_{\tilde D},n_D\}$)
in order to get that the pull-back $\mathcal D_{\tilde Y}$ of $\mathcal D$ to
$\tilde Y$ of [Va2, Subsect. 3.5] is a constant $p$-divisible group isomorphic to
$D\times_k \tilde Y$. The $p$-divisible group $\mathcal D_{\tilde Y}$ is isomorphic
to $\mathcal C/\mathcal L$, where the $p$-divisible group $\mathcal C:=\tilde{\mathcal D}_{m,m}\times_{\mathfrak{i}_D(m)} \tilde Y$ is equipped with an isomorphism
$\theta_1:\mathcal C\to \tilde D\times_k \tilde Y$ and where $\mathcal L$ is the
pull-back of $(D[p^{\kappa}]/\mathcal K)\times_k \tilde Y\subset \tilde
D[p^{\kappa}]$ via an isomorphism $\theta_2:\mathcal C[p^{m-\kappa}]\to \tilde
D[p^{m-\kappa}]\times_k \tilde Y$. Let
$\omega:=\theta_1[p^{m-\kappa}]\circ\theta_2^{-1}\in\pmb{Aut}(\tilde
D[p^{m-\kappa}])(\tilde Y)$. We can assume that the restriction of $\omega$ to the
closed point $\Spec k$ of $\tilde Y$ is the identity. The reduced group $\mathfrak{E}$ of the image of the restriction homomorphism $\pmb{Aut}(\tilde D[p^{m-\kappa}])\to\pmb{Aut}(\tilde D[p^{\kappa}])$ is a finite, \'etale group over $k$ (cf. either Corollary \ref{cor2} (a) for $m\ge 2\kappa+n_{\tilde D}$ or [Va1, Thm. 5.1.1 (c)] for $m>>0$). As the $k$-scheme $\tilde Y$ is connected and reduced (being integral), the restriction of $\omega$ to an automorphism of $\tilde D[p^{\kappa}]_{\tilde Y}$ is a $\tilde Y$-valued point of $\mathfrak{E}$ which is the identity as its restriction to the closed point $\Spec k$ of $\tilde Y$ is so. From this we easily get that the $p$-divisible group $\mathcal D_{\tilde Y}$ is isomorphic to $D\times_k \tilde Y$. 

\section{Generalizations to relative contexts}\label{Sect6}
We recall that $(M,\phi,\vartheta)$ is the contravarint Dieudonn\'e module of $D$ and that we denote also by $\phi:\End(M)[{1\over p}]\to \End(M)[{1\over p}]$ the $\sigma$-linear automorphism induced naturally by $\phi$. Let $G$ be a smooth closed subgroup scheme of $\GL_M$ such that its generic fibre $G_{B(k)}$ is connected. Thus the scheme $G$ is integral. Let $\mathfrak{g}:=\Lie(G)$ be the Lie algebra of $G$. Until the end we will assume that the following two axioms hold for the triple $(M,\phi,G)$:

\medskip
{\bf (AX1)} the Lie subalgebra $\mathfrak{g}[{1\over p}]$ of $\End(M)[{1\over p}]$ is stable under $\phi$, i.e. we have $\phi(\mathfrak{g}[{1\over p}])=\mathfrak{g}[{1\over p}]$;

\smallskip
{\bf (AX2)} there exist a direct sum decomposition $M=F^1\oplus F^0$ such that the following two properties hold:

\begin{definitionlist}
\item the kernel $\bar F^1$ of the reduction modulo $p$ of $\phi$ is $F^1/pF^1$;

\item the cocharacter $\mu:\Bbb G_m\to\GL_M$ which acts trivially on $F^0$ and via the inverse of the identical character of $\Bbb G_m$ on $F^1$, normalizes $G$.
\end{definitionlist}

The triple $(M,\phi,G)$ is called an {\it $F$-crystal with a group} over $k$, cf. [Va1, Def. 1.1 (a) and Subsect. 2.1]. If $\mu$ factors through $G$, then the {\it $W$-condition} of [Va1, Subsubsect. 2.2.1 (d)] holds for $(M,\phi,G)$. Axioms (AX1) and (AX2) are a weakening of the two axioms used in [Va2, Subsect. 4.1]. Until the end, we will use the notations of Subsubsection \ref{dil} and Subsection \ref{act}.

\begin{df}\label{df1}
{\bf (a)} Let $g_1,g_2\in G(W(k))$. An {\it inner isomorphism} between the two triples $(M,g_1\phi,G)$ and  $(M,g_2\phi,G)$ is an element $g_3\in G(W(k))$ that defines an isomorphism between the Dieudonn\'e modules  $(M,g_1\phi,\vartheta g_1^{-1})$ and  $(M,g_2\phi,\vartheta g_2^{-1})$ (i.e. we have $g_3g_1\phi=g_2\phi g_3$). 

{\bf (b)} An {\it inner isomorphism} between $(M/p^mM,g_1[m]\phi_m,\vartheta_m g_1[m]^{-1},G_{W_m(k)})$ and  $(M/p^mM,g_2[m]\phi_m,\vartheta_m g_2[m]^{-1},G_{W_m(k)})$ is an element $g_3[m]\in G(W_m(k))$ which defines an isomorphism between the following two Dieudonn\'e modules $(M/p^mM,g_1[m]\phi_m,\vartheta_m g_1[m]^{-1})$ and  $(M/p^mM,g_2[m]\phi_m,\vartheta_m g_2[m]^{-1})$.
\end{df}

\begin{lemma}\label{L5}
The intersection $\mathcal W^G:=G_k\cap \mathcal W$ is a smooth group over $k$.
\end{lemma} 

\noindent
{\bf Proof:} 
This is a particular case of [CGP, Prop. 2.1.8 (3) and Rm. 2.1.11] applied to the smooth affine group $G_k$ over $k$ and the $\Bbb G_m$ action on it induced by the special fibre $\lambda$ of the inverse $\mu^{-1}$ of the cocharacter $\mu$ of the axiom (AX2) (our pair $(G_k,\mathcal W^G)$ is $(G,P_G(\lambda)=Z_G(\lambda)\times_k U_G(\lambda))$ of loc. cit.).\endproof

\subsection{Relative notations and basic properties} Let $\tilde{\mathcal H}^G$ be the dilatation of $G$ centered on ${\mathcal W}^G$. As in Subsubsection \ref{dil}, based on the axiom (AX1) and Lemma \ref{L5} we get that $\tilde{\mathcal H}^G$ is a smooth affine group scheme over $\Spec W(k)$ equipped with a homomorphism $\tilde{\mathcal P}^G:\tilde{\mathcal H}^G\to G$ whose generic fibre is an isomorphism of groups over $\Spec B(k)$. Moreover, $\tilde{\mathcal H}^G$ is a closed subgroup scheme of $\tilde{\mathcal H}$ (cf. [BLR, Ch. 3, Sect. 3.2, Prop. 2 (a)]) and $\tilde{\mathcal P}^G$ is the natural restriction of $\tilde{\mathcal P}$. Another way to define $\tilde{\mathcal P}^G$ is as follows (cf. Fact \ref{F3} of Appendix A). We can identify $\tilde{\mathcal H}^G$ with the schematic closure $G_{B(k),M,\phi^{-1}(M)}$ in $G\times_{W(k)} \phi^{-1}G\phi$ of $G_{B(k)}$ embedded diagonally in  $G\times_{W(k)} \phi^{-1}G\phi$ (here $\phi^{-1}G\phi$ is the schematic closure of $G_{B(k)}$ in $\GL_{\phi^{-1}(M)}$). Under this identification, $\tilde{\mathcal P}^G$ gets identified with the first projection homomorphism $G_{B(k),M,\phi^{-1}(M)}\to G$.

Let $d_G$ be the relative dimension of $G$ over $\Spec W(k)$; thus $d_G=\dim(G_k)$. Let $m\in\Bbb N^{\ast}$. Let $\mathcal H_m^G:=\Bbb W_m(\tilde{\mathcal H}^G)$; it is a smooth subgroup of $\mathcal H_m$ of dimension $md_G$. Let $\mathcal D_m^G:=\Bbb
W_m(G)$; it is a smooth closed subvariety of $\mathcal D_m$ of dimension $md_G$. 
\begin{df}\label{df2}
{\bf (a)} Let $\pmb{\text{Aut}}(D[p^m])_{\text{crys}}$ and $\zeta_m:\pmb{\text{Aut}}(D[p^m])\to \pmb{\text{Aut}}(D[p^m])_{\text{crys}}$ be as in the proof of Theorem \ref{T5}. Let $\pmb{\text{Aut}}(D[p^m])^G_{\text{crys}}:={\mathcal D}_m^G\cap \pmb{\text{Aut}}(D[p^m])_{\text{crys}}$ (intersection taken inside $\mathcal D_m$). Let $\pmb{\text{Aut}}(D[p^m])^G:=\zeta_{m}^{-1}(\pmb{\text{Aut}}(D[p^m])^G_{\text{crys}})$.
Thus $\pmb{\text{Aut}}(D[p^m])^G(k)$ is the subgroup of $\pmb{\text{Aut}}(D[p^m])(k)$ formed by those elements that define (via $\zeta_{m}$) elements of $G(W_m(k))$. 

{\bf (b)} For $i\in\Bbb N^{\ast}$ let $\gamma_D^G(i):=\dim(\pmb{\text{Aut}}(D[p^i])^G)=\dim(\pmb{\text{Aut}}(D[p^i])^G_{\text{crys}})$. Let $\gamma_D^G(0):=0$. We call $(\gamma_D^G(i))_{i\in\Bbb N}$ the {\it centralizing $G$-sequence} of $D$.
\end{df}

As in Subsection \ref{act} we argue that we have an action 
$$\Bbb T_m^G:\mathcal H_m^G\times_k \mathcal D_m^G\to \mathcal D_m^G$$
which is the natural restriction of the action $\Bbb T_m$. Thus, if $\tilde h\in\tilde{\mathcal H}^G(W(k))$
and $g\in G(W(k))$, then the product of $\tilde h[m]\in \mathcal
H^G_m(k)=\tilde{\mathcal H}^G(W_m(k))$ and $g[m]\in\mathcal D^G_m(k)=G(W_m(k))$
is $\Bbb T_m(\tilde h[m],g[m]):=(hg\phi(h^{-1}))[m]$,
where $h:=\tilde{\mathcal P}^G(W(k))(\tilde h)\in G(W(k))$ has the same meaning as in Subsubsection \ref{dil}. 

Let $n_D^G$ be the smallest non-negative integer that has the following property: for each element $\tilde g\in G(W(k))$ congruent to $1_M$ modulo $p^{n_D^G}$, there exists an inner isomorphism between $(M,\phi,G)$ and  $(M,\tilde g\phi,G)$. The existence of $n_D^G$ is implied by [Va1, Main Thm. A]. If $G=\GL_M$, then it is easy to check based on Lemmas \ref{L1} and \ref{L2} and Corollary \ref{cor4} (a) and (b) that we have $n_D^G=n_D$ (cf. also either [Va1, Lemma 3.2.2 and Cor. 3.2.3] or [NV1, Thm. 2.2 (a)]).  

Let $\rho_{m+1}^G:=\text{Red}_{m+1,\tilde{\mathcal H}^G}$ and
$\tau_{m+1}^G:=\text{Red}_{m+1,G}$ (see Subsubsection \ref{FunctorW_m} for notations).
Thus $\rho_{m+1}^G:\mathcal H_{m+1}^G\twoheadrightarrow\mathcal H_m^G$ is the natural
reduction epimorphism of smooth affine groups over $k$ and $\tau_{m+1}^G:\mathcal D_{m+1}^G\twoheadrightarrow\mathcal D_m^G$ is the natural reduction faithfully flat morphism (epimorphism) of affine (group) varieties over $k$. The action $\Bbb T_m^G$ is also a natural reduction of the action $\Bbb T_{m+1}^G$. The projective
limits 
$$\mathcal H_{\infty}^G:=\text{proj.}\text{lim.}_{m\to\infty} \mathcal
H_m^G\;\;\;\;\;\;\text{and}\;\;\;\;\;\;\mathcal D_{\infty}^G:=\text{proj.}\text{lim.}_{m\to\infty}
\mathcal D_m^G$$ 
in the category of ringed spaces  are also projective
limits in the category of $k$-schemes, $\mathcal H_{\infty}^G$ is an affine group
scheme over $k$, and $\mathcal D_{\infty}^G$ is an affine (group) scheme over $k$. Moreover, we have a limit action
$$\Bbb T_{\infty}^G:\mathcal H_{\infty}^G\times_k\mathcal D_{\infty}^G\to\mathcal
D_{\infty}^G.$$

Let $\mathcal O_m^G$ be the orbit of $1_M[m]\in\mathcal
D_m^G(k)$ under the action $\Bbb T_m^G$. Let $\mathcal S_m^G$ be the subgroup scheme of $\mathcal H_m^G$ which is the stabilizer
of $1_M[m]$. Inside $\mathcal H_m$, we have $\mathcal S_m^G=\mathcal S_m\cap \mathcal H_m^G$. Let $\mathcal C_m^G:=\mathcal S_{m,\re}^G$. We have $\dim(\mathcal
S_m^G)=\dim(\mathcal C_m^G)=\dim(\mathcal H_m^G)-\dim(\mathcal O_m^G)=md_G-\dim(\mathcal O_m^G)$. The orbits of $\Bbb T_{\infty}^G$ are defined similarly to the orbits of $\Bbb T_{\infty}$ (see Subsubsection \ref{limit}).

\begin{lemma}\label{L6}
The set of orbits of $\Bbb T_m^G$ is in natural bijection to the set of (representatives of) inner isomorphism classes of quadruples of the form $\break (M/p^mM,g[m]\phi_m,\vartheta_m g[m]^{-1},G_{W_m(k)})$  with $g\in G(W(k))$. 
\end{lemma}

\noindent
{\bf Proof:} The rule that maps the orbit of $g[m]$ under $\Bbb T_m^G$ to the inner isomorphism class of $(M/p^mM,g[m]\phi_m,\vartheta_m g[m]^{-1},G_{W_m(k)})$ is well defined (cf. definition of $\Bbb T_m^G$) and we are left to check that the resulting map is a bijection. It suffices to show that if $g_1,g_2,g_3\in G(W(k))$ are as in Definition \ref{df1} (b), then there exists $\tilde h\in \tilde{\mathcal H}^G(W(k))$ such that $\Bbb T_m^G(\tilde h[m],g_1[m])=g_2[m]$. The pair $(g_3[m],(g_2^{-1}g_3g_1)[m])$ is an automorphism of the pair of maps 
$\xymatrixcolsep{0.8pc}\xymatrix@1{ M/p^mM \ar@<0.5ex>[rr]^{\phi_m} & &M/p^mM
\ar@<0.5ex>[ll]^{\vartheta_m} \\}$. Thus the pair $(g_3[m],\phi^{-1}(g_2^{-1}g_3g_1)[m]\phi)$ is an isomorphism of the pair of $W_m(k)$-linear maps 
$\xymatrixcolsep{0.8pc}\xymatrix@1{ M/p^mM \ar@<0.5ex>[rr]^{j_m} & &\phi^{-1}(M)\ar@<0.5ex>[ll]^{p} \\}$$/p^m\phi^{-1}(M)$, where $j_m$ is the reduction modulo $p^m$ of the inclusion $M\hookrightarrow \phi^{-1}(M)$. 

Based on the axiom (AX2), the group $G_{B(k)}$ is well positioned with respect to $(M,\phi^{-1}(M))$ in the sense of Definition \ref{df3} of Appendix A (cf. Proposition \ref{P4} of Appendix A applied to $(W(k),p,G,\mu_{B(k)},\phi^{-1}(M))$ instead of $(R,\pi,G_M,\lambda,M')$). For the remaining part of the proof we will only use this well positioned property which asserts the existence of an isomorphism $\Theta:\tilde{\mathcal H}^G=G_{B(k),M,\phi^{-1}(M)}\to G_{B(k),M,\phi^{-1}(M)}^*:=(G\times_{W(k)} \phi^{-1}G\phi)\cap \GL_{M,\phi^{-1}(M)}$.

We have $(g_3[m],\phi^{-1}(g_2^{-1}g_3g_1)[m]\phi)\in G_{B(k),M,\phi^{-1}(M)}^*(W_m(k))$, cf. first before last paragraph and axiom (AX1) which implies that $\phi^{-1}g_2^{-1}g_3g_1\phi\in G(B(k))$. Let $\tilde h\in\tilde{\mathcal H}^G(W(k))$ be an element that lifts the unique element $\tilde h[m]\in\tilde{\mathcal H}^G(W_m(k))$ for which we have $\Theta(\tilde h[m])=(g_3[m],\phi^{-1}(g_2^{-1}g_3g_1)[m]\phi)$ (cf. previous paragraph). Therefore $\tilde{\mathcal P}^G(\tilde h[m])=g_3[m]$ and the linear automorphisms of $\phi^{-1}(M)/p^m\phi^{-1}(M)$ induced by $\tilde h$ and $\phi^{-1}(g_2^{-1}g_3g_1)[m]\phi$ coincide (i.e. $\phi(\tilde h)[m]=(g_2^{-1}g_3g_1)[m]$). Thus we compute $\Bbb T_m^G(\tilde h[m],g_1[m])=g_3[m]g_1[m]((g_2^{-1}g_3g_1)[m])^{-1}=g_2[m]$.\endproof

\medskip
There exist upper bounds of $n_D^G$ which depend only on $r$ and on the closed embedding homomorphism $G\hookrightarrow \GL_M$, cf. [Va1, Main Thm. A and Example 3.1.5]. The following consequences of Lemma \ref{L6} and of the existence of such upper bounds of $n_D^G$ are shown similarly to Corollary \ref{cor4} (a) to (c). 

\begin{cor}\label{cor7} The following three properties hold:

\medskip
{\bf (a)} The projective limit $\mathcal O_{\infty}^G:=\text{proj.}\text{lim.}_{m\to\infty} \mathcal O_m^G$ in the category of ringed spaces is also a projective limit in the category of $k$-schemes and it is a reduced, quasi-compact, locally closed subscheme of $\mathcal D_{\infty}^G$ whose $k$-valued points form one orbit of the action $\Bbb T_{\infty}^G(k):\mathcal H_{\infty}^G(k)\times\mathcal D_{\infty}^G(k)\to \mathcal D_{\infty}^G(k)$ in the category of sets. Thus $\mathcal O_{\infty}^G$ is the orbit of $1_M\in\mathcal D_{\infty}^G(k)$ under the action $\Bbb T_{\infty}^G$.

 \smallskip
{\bf (b)} The set of orbits of $\Bbb T_{\infty}^G$ is in natural bijection to the set of (representatives of) inner isomorphism classes of triples  of the form $(M,g\phi,\vartheta g,G)$  with $g\in G(W(k))$.

\smallskip
{\bf (c)} If $K$ is an algebraically closed field that contains $k$, then we have $n_{D_K}^{G_{W(K)}}=n_D^G$.
\end{cor}

We note that the proof of Corollary \ref{cor7} (b) does not require the use of the well positioned property mentioned in the proof of Lemma \ref{L6}.

The generalizations of Theorems \ref{UP}, \ref{CP}, and \ref{INF} to the relative context provided by $(M,\phi,G)$ is automatic. For instance, if $\mathcal{T}_{m+1}^G$ is the reduced group of $(\rho_{m+1}^G)^{-1}(\mathcal{S}_m^G)$ and if $\mathcal{V}_{m+1}^G$ is the inverse image of the point $1_M[m]\in\mathcal{O}_m^G(k)$ under the epimorphism $\tau_{m+1}^G:\mathcal{D}_{m+1}^G\twoheadrightarrow\mathcal{D}_m^G$, then the set of orbits of the action $\Bbb T_{m+1}^G$ that map onto $\mathcal{O}_m^G$ is in a natural bijection to the set of orbits of the action $\mathcal{T}_{m+1}^G\times_k \mathcal{V}_{m+1}^G\to \mathcal{V}_{m+1}^G$ induced via restriction by $\Bbb T_{m+1}^G$. Moreover, $\mathcal{V}_{m+1}^G$ is isomorphic to $\Bbb A_k^{d_G}$, $\mathcal{T}_{m+1}^{G,0}$ is unipotent, we have a short exact sequence $1\to \Bbb G_a^{d_G}\to\mathcal{T}_{m+1}^G\to\mathcal{C}_{m+1}^G\to 1$, etc. Thus in all that follows we will concentrate on the generalization of Theorem \ref{Gamma} and of its corollaries.

\begin{lemma}\label{L7}
We assume that $n_D=0$ (i.e. $D$ is either \'etale or of multiplicative type). Then we have $n_D^G\in\{0,1\}$. If moreover the special fibre $G_k$ is connected, then we have $n_D^G=0$. 
\end{lemma}

\noindent
{\bf Proof:}
As $n_D=0$, the image of $\mu$ is either trivial or equal to the center of $\GL_M$. Thus $\mathcal{W}=\GL_{M/pM}$, $\mu$ centralizes $G$, and we have $\mathcal{W}^G=G_k$. Let $\sigma_{\phi}:=\phi\mu(p):M\to M$; it is a $\sigma$-linear automorphism of $M$ which normalizes $G$ and which acts on $G(B(k))$ in the same way as $\phi$. The $\Bbb Z_p$-module $M_{\Bbb Z_p}:=\{x\in M|\sigma_{\phi}(x)=x\}$ is a $\Bbb Z_p$ structure of $M$ (i.e. we have $M=W(k)\otimes_{\Bbb Z_p} M_{\Bbb Z_p}$) and $\mathfrak{g}_{\Bbb Z_p}:=\{x\in\mathfrak{g}|\sigma_{\phi}(x)=x\}=\mathfrak{g}\cap\End(M_{\Bbb Z_p})$ is a $\Bbb Z_p$ structure of $\mathfrak{g}$. Thus $G$ is the pull-back to $\Spec W(k)$ of a smooth closed subgroup scheme $G_{\Bbb Z_p}$ of $\GL_{M_{\Bbb Z_p}}$, cf. [Va4, Prop. 3.2]. 

If $G_k$ is (resp. is not) connected, let $g\in G(W(k))$ (resp. let $g\in\Ker(G(W(k))\to G(k))$). We will show that there exist an inner isomorphism between $(M,g\phi,G)$ and $(M,\phi,G)$. For this let $s$ be a positive integer such that $s\ge\max\{2,n_D^G\}$. The affine group $\Bbb K_s:=\Bbb W_s(G)$ (resp. the kernel $\Bbb K_s$ of the natural restriction homomorphism $\Bbb W_s(G)\to \Bbb W_1(G)$) over $k$ is connected and smooth, cf. Subsubsection \ref{FunctorW_m}. The endomorphism of $\Bbb K_s$ induced naturally by $\sigma_{\phi}$ is the usual Frobenius endomorphism $\sigma_s$ of $\Bbb K_s$ with respect to the $\Bbb F_p$-form $\Bbb K_{s,\Bbb F_p}$ of $\Bbb K_s$ defined by $\Bbb W_s(G_{\Bbb Z_p})$ (resp. by the kernel of the natural restriction homomorphism $\Bbb W_s(G_{\Bbb Z_p})\to \Bbb W_1(G_{\Bbb Z_p})$); here the $\Bbb W_s$ and $\Bbb W_1$ functors are applied over $\Bbb Z_p$. A classical theorem of Lang implies that there exists an element $h[s]\in \Bbb K_s(k)\subset G(W_s(k))$ such that $g[s]=h[s]^{-1}\sigma_{s}(h[s])$. If $h\in G(W(k))$ lifts $h[s]$, then we get that $g_1:=hg\phi(h^{-1})\in G(W(k))$ is such that $g_1[s]=1_M[s]$ and $hg\phi=hg\phi(h^{-1})\phi h=g_1\phi h$. From this and the definition of $n_D^G\le s$, we get that there exists $h_1\in G(W(k))$ such that $h_1g_1\phi=\phi h_1$. We conclude that $h_1h$ is an inner isomorphism between $(M,g\phi,G)$ and $(M,\phi,G)$.
Thus we have $n_D^G=0$ (resp. $n_D^G\le 1$). \endproof

\begin{prop}\label{P2}
The following five properties hold:

\medskip
{\bf (a)} The natural homomorphism $\iota_m^G:\mathcal S_m^G\to \pmb{Aut}(D[p^m])^G_{\text{crys}}$ induced by $\iota_m$ gives an isomorphism $\iota_{m}^G(k):\mathcal S_m^G(k)\to \pmb{Aut}(D[p^m])^G_{\text{crys}}(k)$.

\smallskip
{\bf (b)} The sequence $(\gamma_D^G(i))_{i\in\Bbb N}$ is increasing.

\smallskip
{\bf (c)} The subsequence $(\gamma_D^G(i))_{i\ge n_D^G}$ is constant. 

\smallskip
{\bf (d)} We assume that $n_D^G\ge 2$. Then we have $\gamma_D^G(n_D^G-1)<\gamma_D^G(n_D^G)$. 

\smallskip
{\bf (e)} The connected smooth affine groups $\pmb{Aut}(D[p^m])^{G,0}_{\re}$ and $\pmb{Aut}(D[p^m])^{G,0}_{\text{crys},\re}$ are unipotent. 
\end{prop}

\noindent
{\bf Proof:}  Due to the axiom (AX2), $\sigma_{\phi}:=\phi\mu(p):M\to M$ is a $\sigma$-linear automorphism of $M$ which normalizes $G$. Based on this, the argument for (a) is the same as the one in the proof of [Va2, Cor. 4.3 (b)] (see second paragraph of [Va2, p. 635]). As a second proof of (a),  we first remark that $\iota_m^G(k)$ is injective as $\iota_m(k)$ is so; the surjectivity of $\iota_m^G(k)$ follows directly from the proof of Lemma \ref{L6} applied with $g_1[m]=g_2[m]=1_M[m]$.

Based on (a) we get: 
$$\gamma_D^G(m)=\dim(\pmb{Aut}(D[p^m])^G_{\text{crys}})=\dim(\mathcal C_m^G)=md_G-\dim(\mathcal O_m^G).$$ 
Thus $\gamma_D^G(m+1)-\gamma_D^G(m)=d_G-\dim(\mathcal O_{m+1}^G)+\dim(\mathcal O_m^G)$. As the fibres of $\tau_{m+1}^G$ are isomorphic to $\Bbb A_k^{d_G}$, the faithfully flat morphism $\mathcal O_{m+1}^G\to \mathcal O_m^G$ induced by $\tau_{m+1}^G$ has fibres of dimension at most $d_G$. Thus $\dim(\mathcal O_{m+1}^G)-\dim(\mathcal O_m^G)\le d_G$, i.e. $\gamma_D^G(m+1)-\gamma_D^G(m)\ge 0$. As $\gamma_D^G(1)\ge 0=\gamma_D^G(0)$, we get that (b) holds.

If $n_D^G=0$, then $\mathcal{O}_m^G=\mathcal{D}_m^G$ has dimension $md_G$ and thus $\gamma_D^G(m)=0$. As $\gamma_D^G(0)=0$, we conclude that (c) holds if $n_D^G=0$. Thus to prove (c) we can assume that $n_D^G\ge 1$. If $m\ge n_D^G$, then from the very definition of $n_D^G$ we get that $\mathcal O_{m+1}^G=(\tau_{m+1}^G)^{-1}(\mathcal O_m^G)$. This implies that $\dim(\mathcal O_{m+1}^G)=\dim(\mathcal O_m^G)+d_G$, i.e. $\gamma_D^G(m+1)-\gamma_D^G(m)=0$. Thus (c) holds. 

To prove (d), it suffices to show that the equality $\gamma_D^G(n_D^G-1)=\gamma_D^G(n_D^G)$ leads to a contradiction. This equality implies that $\dim(\mathcal O_{n_D^G}^G)=d_G+\dim(\mathcal O_{n_D^G-1}^G)$. Therefore all fibres of the faithfully flat morphism $\mathcal O_{n_D^G}^G\to\mathcal O_{n_D^G-1}^G$ have dimension $d_G$ and thus, being closed subschemes of $\Bbb A_k^{d_G}$, are isomorphic to $\Bbb A_k^{d_G}$. Thus $\mathcal O_{n_D^G}^G=(\tau_{n_D^G}^G)^{-1}(\mathcal O_{n_D^G-1}^G)$. This means that if $\tilde g\in\Ker(G(W(k))\to G(W_{n_D^G-1}(k)))$, then $(M,\tilde g\phi,G)$ is inner isomorphic to a triple of the form $(M,\tilde g_1\phi,G)$ with $\tilde g_1\in\Ker(G(W(k))\to G(W_{n_D^G}(k)))$ and therefore (cf. definition of $n_D^G$) it is inner isomorphic to $(M,\phi,G)$. Thus we have $n_D^G\le n_D^G-1$. Contradiction. Therefore (d) holds. 

Part (e) follows from the fact that $\pmb{Aut}(D[p^m])^{G,0}_{\re}$ and $\pmb{Aut}(D[p^m])^{G,0}_{\text{crys},\re}$ are subgroups of the unipotent groups $\pmb{Aut}(D[p^m])^{0}_{\re}$ and $\pmb{Aut}(D[p^m])^{0}_{\text{crys},\re}$ (respectively), cf. Theorem \ref{T5} (b) and (c) and the definitions. 
\endproof

\subsection{Lie algebra group schemes of endomorphisms}\label{Lie}
For all that follows, it is convenient to also introduce the Lie algebra version of $\pmb{Aut}(D[p^m])^G_{\text{crys}}$. Let $l,i\in\Bbb N$ be such that $l\ge
i$. Let $\underline{\sharp}$ be the vector group scheme over $\Spec W(k)$ (resp. $\Spec W_m(k)$) associated to a free $W(k)$-module (resp. $W_m(k)$-module) $\sharp$ of finite rank. Let $\Bbb W_0(\underline{\mathfrak{g}}):=\Spec k$. 

We recall that $\pmb{End}(D[p^l])_{\text{crys}}$ is a closed subgroup scheme of $\Bbb W_l(\underline{\End(M)})$, cf. its definition before Theorem \ref{T5}. Let  $\pmb{End}(D[p^l])_{\text{crys}}^{\mathfrak{g}}:=\pmb{End}(D[p^l])_{\text{crys}}\cap \Bbb W_l(\underline{\mathfrak{g}})$, the intersection being taken inside $\Bbb W_l(\underline{\End(M)})$. The crystalline Dieudonn\'e functor defines a homomorphism
$\Bbb D:\pmb{End}(D[p^l])\to \pmb{End}(D[p^l])_{\text{crys}}$
which induces a bijection at the level of $k$-valued points. Let $\pmb{End}(D[p^l])^{\mathfrak{g}}:=\Bbb D^{-1}(\pmb{End}(D[p^l])_{\text{crys}}^{\mathfrak{g}})$. The $k$-valued points of $\pmb{End}(D[p^l])^{\mathfrak{g}}$ are those endomorphisms of $D[p^l]$ whose crystalline realizations belong to $\mathfrak{g}/p^l\mathfrak{g}\subset \End(M/p^lM)$. 

The homomorphism $\pmb{End}(D[p^l])^{\mathfrak{g}}\to \pmb{End}(D[p^l])_{\text{crys}}^{\mathfrak{g}}$ which is the natural restriction of $\Bbb D$, induces a bijection at the level of $k$-valued points. Thus let 
$$\gamma_D^{\mathfrak{g}}(l):=\dim(\pmb{End}(D[p^l])^{\mathfrak{g}})=\dim(\pmb{End}(D[p^l])_{\text{crys}}^{\mathfrak{g}}).$$
\noindent
We have $\gamma_D^{\mathfrak{g}}(0)=0$. 

Next we will use the notations of Subsection \ref{SSect2.1}. Let 
$$\kappa_{i,l}^{\mathfrak{g}}:\pmb{End}(D[p^i])^{\mathfrak{g}}\hookrightarrow \pmb{End}(D[p^l])^{\mathfrak{g}}$$ 
be the monomorphism induced naturally by $\kappa_{i,l}$. Let 
$$r_{l,i}^{\mathfrak{g}}:\pmb{End}(D[p^l])^{\mathfrak{g}}\to \pmb{End}(D[p^i])^{\mathfrak{g}}$$ 
be the restriction homomorphism induced naturally by $r_{l,i}$. As in Subsection \ref{SSect2.1} we argue that we have a complex of commutative group schemes
 $$0\to \pmb{End}(D[p^i])^{\mathfrak{g}}\xrightarrow{\kappa_{i,l}^{\mathfrak{g}}}
\pmb{End}(D[p^l])^{\mathfrak{g}}\xrightarrow{r_{l,l-i}^{\mathfrak{g}}} \pmb{End}(D[p^{l-i}])^{\mathfrak{g}}$$
whose complex of groups of $k$-valued points is exact. This implies that:

\medskip
{\bf (i)} The image $\text{Im}(r_{l,l-i}^{\mathfrak{g}})$ has dimension $\gamma_D^{\mathfrak{g}}(l)-\gamma_D^{\mathfrak{g}}(i)$.

\medskip
If $l\ge n_D$, then $\text{Im}(r_{l+1,1})$ and thus also $\text{Im}(r_{l+1,1}^{\mathfrak{g}})$ has dimension $0$ (cf. Corollary \ref{cor2} (a)). From this and (i) we get that:

\medskip
{\bf (ii)} If $l\ge n_D$, then $\gamma_D^{\mathfrak{g}}(l+1)=\gamma_D^{\mathfrak{g}}(l)$.

\begin{thm}\label{T6}
We assume that one of the following two conditions holds:

\medskip
{\bf (i)} The subset $\mathfrak{g}$ of $\End(M)$ is stable under products (therefore $W(k)1_M+\mathfrak{g}$ is a $W(k)$-subalgebra of $\End(M)$ and $G$ is the group scheme of invertible elements of $1_M+\mathfrak{g}$).

\smallskip
{\bf (ii)} We have $p>2$ and $G$ admits a Cayley transform $\mathcal{U}$ of the form $g\mapsto (1_M-g)(1_M+g)^{-1}$ (i.e. the rule $g\mapsto (1_M-g)(1_M+g)^{-1}$ defines an isomorphism $\mathcal{U}$ between the non-empty open subscheme of $G$ formed by points which do not have the eigenvalue $-1$ and the non-empty open subscheme of $\underline{\mathfrak{g}}$ formed by points which do not have the eigenvalue $-1$). If $n_D=0$, then moreover $G_k$ is connected. 

\medskip
Then the following seven properties hold:

\medskip
{\bf (a)} For all $l\in\Bbb N$ we have $\gamma_D^G(l)=\gamma_D^{\mathfrak{g}}(l)$.

\smallskip
{\bf (b)} The sequence $(\gamma^G_D(i+1)-\gamma^G_D(i))_{i\in\Bbb N}$ is a decreasing sequence in $\Bbb N$.

\smallskip
{\bf (c)} If $n_D^G>0$, then we have $\gamma_D^G(1)<\gamma_D^G(2)<\cdots<\gamma_D^G(n_D^G)$.

\smallskip
{\bf (d)} We have $n_D^G\le n_D$.

\smallskip
{\bf (e)} Let $l>i>0$ be integers. Then the image of the restriction homomorphism $r_{l,i}^{\mathfrak{g}}:\pmb{End}(D[p^l])^{\mathfrak{g}}\to \pmb{End}(D[p^i])^{\mathfrak{g}}$ (or $r_{l,i,\text{crys}}^{\mathfrak{g}}:\pmb{End}(D[p^l])_{\text{crys}}^{\mathfrak{g}}\to \pmb{End}(D[p^i])_{\text{crys}}^{\mathfrak{g}}$) is finite if and only if $l-i\ge n_D^G$.

\smallskip
{\bf (f)} Let $l$ and $i$ be as in (e). Then the image of the restriction homomorphism $s_{l,i}^G:\pmb{Aut}(D[p^l])^G\to \pmb{Aut}(D[p^i])^G$ (or $s_{l,i,\text{crys}}^G:\pmb{Aut}(D[p^l])_{\text{crys}}^G\to \pmb{Aut}(D[p^i])_{\text{crys}}^G$) is finite if and only if $l-i\ge n_D^G$.

\smallskip
{\bf (g)} If $i\ge 2n_D^G$ is a positive integer, then the unipotent groups $\pmb{Aut}(D[p^i])_{\re}^{G,0}$ and $\pmb{Aut}(D[p^i])_{\text{crys},\re}^{G,0}$ are commutative. 
\end{thm}

\noindent
{\bf Proof:} 
To prove (a) we distinguish three cases as follows.

We assume that (i) holds and $1_M\in\mathfrak{g}$ (i.e. $\mathfrak{g}$ is a $W(k)$-subalgebra of $\End(M)$). In this case $\pmb{Aut}(D[p^m])_{\text{crys}}^G$ is an open subscheme of $\pmb{End}(D[p^m])_{\text{crys}}^{\mathfrak{g}}$ and thus (a) holds. 

We assume that (i) holds and $1_M\notin\mathfrak{g}$. Then $\pmb{Aut}(D[p^m])^G_{\text{crys}}$ is isomorphic to an open subscheme of $\pmb{End}(D[p^m])_{\text{crys}}^{\mathfrak{g}}$ via the morphism which at the level of $k$-valued points maps an automorphism $g[m]$ to the endomorphism $g[m]-1_M[m]$. Thus (a) holds.

We assume that (ii) holds. To check that (a) holds, we can assume $l>0$. The isomorphism $\mathcal{U}$ induces naturally a Cayley transform $\mathcal{U}_l:OG_l\to O\mathfrak{g}_l$ which is an isomorphism between  the non-empty open subscheme $OG_l$ of $G_{W_l(k)}$ formed by points which do not have the eigenvalue $-1$ and the non-empty open subscheme $O\mathfrak{g}_l$ of $\underline{\mathfrak{g}/p^l\mathfrak{g}}$ formed by points which do not have the eigenvalue $-1$. 

The association $\diamond\mapsto (1-\diamond)(1+\diamond)^{-1}$ maps endomorphisms of some pull-back of $D[p^l]$ or of some scalar extension of $(M/p^lM,\phi_l,\vartheta_l)$  to endomorphisms of the same pull-back of $D[p^l]$ or of the same scalar extension of $(M/p^lM,\phi_l,\vartheta_l)$. Therefore $\mathcal{U}_l$ induces an isomorphism 
$$\mathcal{U}_l:\pmb{OAut}(D[p^l])^G_{\text{crys}}\to\pmb{OEnd}(D[p^l])^{\mathfrak{g}}_{\text{crys}}$$ 
from a suitable non-empty open subscheme $\pmb{OAut}(D[p^l])^G_{\text{crys}}$ of $\pmb{Aut}(D[p^l])^G_{\text{crys}}$ to a suitable non-empty open subscheme $\pmb{OEnd}(D[p^l])^{\mathfrak{g}}_{\text{crys}}$ of $\pmb{End}(D[p^l])^{\mathfrak{g}}_{\text{crys}}$. Thus $\dim(\pmb{Aut}(D[p^l])^G_{\text{crys}})=\dim(\pmb{End}(D[p^l])^{\mathfrak{g}}_{\text{crys}})$. Therefore (a) holds. 

Based on the property \ref{Lie} (i) and on (a), the proofs of  (b) and (c) are the same as the proofs of Theorem \ref{Gamma} (a) and (b) (the roles of $\pmb{End}(D[p^m])$, $\pmb{Aut}(D[p^m])$, and $\Bbb T_m$ being replaced by $\pmb{End}(D[p^m])^{\mathfrak{g}}$, $\pmb{Aut}(D[p^m])^G$, and $\Bbb T_m^G$).

If (i) holds, then $G_k$ is the group scheme of invertible elements of $1_{M/pM}+\mathfrak{g}/p\mathfrak{g}$ and thus it is connected. Thus to prove (d) we can use the fact that if $n_D=0$, then $G_k$ is connected. If $n_D=0$, then $n_D^G=0$ by Lemma \ref{L7}. Thus to prove (d) we can assume that $\min\{n_D,n_D^G\}>0$. Let $l:=n_D+1\geq 2$. We have $\gamma_D^{\mathfrak{g}}(l)=\gamma_D^{\mathfrak{g}}(l-1)$, cf. property \ref{Lie} (ii). From this and (a) we get $\gamma_D^G(l)=\gamma_D^G(l-1)$. From this and (b) we get that $l=n_D+1\ge n_D^G+1$.  Thus (d) holds. 

Based on (a) to (c), part (e) follows from the property \ref{Lie} (i). 

Part (f) follows from (e) via the three cases considered in the proof of (a). For instance, if the condition (ii) holds, then the existence of $\mathcal{U}_l:\pmb{OAut}(D[p^l])^G_{\text{crys}}\to\pmb{OEnd}(D[p^l])^{\mathfrak{g}}_{\text{crys}}$ implies that $r_{l,i,\text{crys}}^{\mathfrak{g}}$ has finite image if and only if $s_{l,i,\text{crys}}^G$ has finite image. Thus indeed (e) implies (f). 

The argument that $\pmb{Aut}(D[p^i])_{\re}^{G,0}$ (and thus also $\pmb{Aut}(D[p^i])_{\text{crys},\re}^{G,0}$) is commutative if $i\ge \max\{1,2n_D^G\}$ is the same as the one for Corollary \ref{cor3} (d). Thus (g) holds.
\endproof

\subsection{The invariant $\upsilon_D$, exponentials, and logarithms}\label{inv}
In order to obtain some variants of Theorem \ref{T6}, we will introduce new invariants $\upsilon_D$ and $\upsilon_D^G$. The connected smooth affine group $\pmb{Aut}(D[p])_{\text{crys},\re}^0$ is unipotent, cf. Theorem \ref{T5} (a) and (b). Thus $\pmb{Aut}(D[p])_{\text{crys},\re}^0$ is a subgroup of the unipotent radical of some Borel subgroup of $\GL_{M/pM}$. This implies that there exists a smallest non-negative integer $\upsilon_D$ with the property that for each $g[1]\in \pmb{Aut}(D[p])_{\text{crys}}^0(k)$, we have $(g[1]-1_M[1])^{\upsilon_D}=0$. We call $\upsilon_D$ as the {\it unipotent number} of $D$. We recall that $r$ is the height of $D$ (i.e. the rank of $M$). We have $\upsilon_D\in\{0,\ldots,r\}$. From Corollary \ref{cor6} (b) we get: 

\medskip
{\bf (i)} For each element $e[1]\in\pmb{End}(D[p])_{\text{crys}}^0(k)$ we have $e[1]^{\upsilon_D}=0$. 

\medskip
Let $\upsilon_D^G$ be the smallest non-negative integer with the property that for each $g[1]\in \pmb{Aut}(D[p])^{G,0}_{\text{crys}}(k)$ and each $e[1]\in \pmb{End}(D[p])^{\mathfrak{g},0}_{\text{crys}}(k)$ we have $(g[1]-1_M[1])^{\upsilon_D^G}=0$ and $e[1]^{\upsilon_D^G}=0$ as elements of $\End(M/pM)$. We call $\upsilon_D^G$ as the {\it unipotent $G$-number} of $D$. Obviously:

\medskip
{\bf (ii)} We have inequalities $0\le \upsilon_D^G\le \upsilon_D\le r$. 

\medskip
For $x\in\End(M)$, let $\bar x\in \End(M/pM)$ be its reduction modulo $p$. For $t\in\{1,\ldots,r\}$ let $\Sigma_{1,t}(M):=p\{x\in\End(M)|\bar x^t=0\}$. For $p>2$ and $s\in\{1,\ldots,\min\{r,{{p-1}\over 2}\}\}$, let $\Sigma_s(M):=\{x\in\End(M)|\bar x^s=0\}$. Let $\Sigma$ denote either $\Sigma_{1,t}$ or $\Sigma_s$. The exponential rule $z\mapsto \sum_{i=0}^{\infty} {{z^i}\over {i!}}$ defines a bijection
$$\exp:\Sigma(M)\to 1_M+\Sigma(M)$$
whose inverse 
$$\log:1_M+\Sigma(M)\to \Sigma(M)$$ 
is defined by the logarithmic rule $1_M+z\mapsto \sum_{i=1}^{\infty} {{(-1)^{i-1}z^i}\over {i}}$, cf. Proposition \ref{P6} (b) or \ref{P7} (a) of Appendix B. Let $\Sigma^{\mathfrak{g}}:=\Sigma(M)\cap\mathfrak{g}$ and $\Sigma^G:=\{g\in G(W(k))|g-1_M\in\Sigma(M)\}$. Let $\Sigma^{\mathfrak{g}}[m]$ (resp. $\Sigma^G[m]$) be the image of $\Sigma^{\mathfrak{g}}$ (resp. of $\Sigma^G$) in $\mathfrak{g}/p^m\mathfrak{g}$ (resp. in $G(W_m(k))$). From Proposition \ref{P7} (a) of Appendix B we get:

\medskip
{\bf (iii)} Restricting $\exp$ to $\Sigma^{\mathfrak{g}}$ and $\log$ to $\Sigma^G$ we get inverse bijections $\xymatrixcolsep{0.6pc}\xymatrix@1{ \Sigma^{\mathfrak{g}} \ar@<0.5ex>[rr]^{\exp} & &\Sigma^G
\ar@<0.5ex>[ll]^{\log} \\}$ which induce inverse bijections 
$\xymatrixcolsep{0.9pc}\xymatrix@1{ \Sigma^{\mathfrak{g}}[m] \ar@<0.5ex>[rr]^{\exp_m} & &\Sigma^G[m]
\ar@<0.5ex>[ll]^{\log_m} \\}$.

\medskip
Let $\mathcal D_{m,\Sigma}^G$ be the reduced closed subgroup scheme of $\mathcal D_{m}^G$ whose $k$-valued points are $\Sigma^G[m]$. Let $\Bbb W_{m,\Sigma}(\underline{\mathfrak{g}})$ be the reduced closed subgroup scheme of $\Bbb W_{m}(\underline{\mathfrak{g}})$ whose $k$-valued points are $\Sigma^{\mathfrak{g}}[m]$. Let $\Omega^G_m:=\Ker(s_{m,1,\text{crys}}^G)_{\re}^0$, where $s_{m,1,\text{crys}}^G:\pmb{Aut}(D[p^m])^G_{\text{crys}}\to \pmb{Aut}(D[p])^G_{\text{crys}}$ is the natural restriction homomorphism; it is a closed subgroup of $\mathcal D_{m,1}^G:=\Ker(\tau_2^G\circ\cdots\circ\tau_m^G:\mathcal D_m^G\to\mathcal D_1^G)$. Let $\Omega^{\mathfrak{g}}_m:=\Ker(r_{m,1,\text{crys}}^{\mathfrak{g}})_{\re}^0$, where $r_{m,1,\text{crys}}^{\mathfrak{g}}:\pmb{End}(D[p^m])^{\mathfrak{g}}_{\text{crys}}\to \pmb{End}(D[p])^{\mathfrak{g}}_{\text{crys}}$ is the natural restriction homomorphism; it is a closed subgroup scheme of $\Bbb W_{m,1}(\underline{\mathfrak{g}}):=\Ker(\Bbb W_m(\underline{\mathfrak{g}})\to \Bbb W_1(\underline{\mathfrak{g}}))$.  

If $e_m\in \Sigma^{\mathfrak{g}}[m]$ (resp. $e_m\in \Sigma^G[m]$) is the crystalline realization of an endomorphism (resp. automorphism) of $D[p^m]$, i.e. is a $k$-valued point of $\pmb{End}(D[p^m])_{\text{crys}}^{\mathfrak{g}}$ (resp. of $\pmb{Aut}(D[p^m])_{\text{crys}}^G$), then $\exp_m(e_m)\in \Sigma^G[m]$ (resp. $\log_m(e_m)\in \Sigma^{\mathfrak{g}}[m]$) is the crystalline realization of an automorphism (resp. endomorphism) of $D[p^m]$ (cf. Proposition \ref{P7} (b) of Appendix B). From this, the property (iii), and the Proposition \ref{P6} (c) of Appendix B we get the following two properties.

\medskip
{\bf (iv)} If $p$ is arbitrary and $\Sigma=\Sigma_{1,t}$, then there exists an isomorphism $\mathcal D_{m,\Sigma}^G\to \Bbb W_{m,\Sigma}(\underline{\mathfrak{g}})$ which on $k$-valued points is induced by $\log_m$ and whose inverse is induced on $k$-valued points by $\exp_{m}$. For $t\ge\upsilon_D$ the isomorphism $\mathcal D_{m,\Sigma}^G\to \Bbb W_{m,\Sigma}(\underline{\mathfrak{g}})$ restricts to an isomorphism $\Omega_m^G\to\Omega_m^{\mathfrak g}$. If $p>2$, then the isomorphism $\mathcal D_{m,\Sigma}^G\to \Bbb W_{m,\Sigma}(\underline{\mathfrak{g}})$ is the restriction of an isomorphism $\mathcal D_{m,1}^G\to \Bbb W_{m,1}(\underline{\mathfrak{g}})$ which on $k$-valued points is induced by $\log_m$ and whose inverse is induced on $k$-valued points by $\exp_{m}$.

\smallskip
{\bf (v)} If $p>2$ and  $\Sigma=\Sigma_s$ with $s\in\{1,\ldots,\min\{r,{{p-1}\over 2}\}\}$, then there exists an isomorphism $\mathcal D_{m,\Sigma}^G\to \Bbb W_{m,\Sigma}(\underline{\mathfrak{g}})$ which on $k$-valued points is induced by $\log_m$ and whose inverse is induced on $k$-valued points by $\exp_{m}$. If $\upsilon_D^G\le {{p-1}\over 2}$, then for $s\ge \upsilon_G^D$ the mentioned isomorphism restricts to an isomorphism $\pmb{Aut}(D[p^m])^{G,0}_{\text{crys},\re}\to \pmb{End}(D[p^m])^{\mathfrak{g},0}_{\text{crys},\re}$.

\begin{prop}\label{P3}
The following two properties hold:

\medskip
{\bf (a)} If $l\ge n_D+1$, then $\gamma_D^G(l)=\gamma_D^{\mathfrak{g}}(l)$.

\smallskip
{\bf (b)} We have $n_D^G\le n_D+1$. 
\end{prop}

\noindent
{\bf Proof:}  If $l\ge n_D+1$, then from Corollary \ref{cor2} (a) we get that we have identities $\Omega_l^G=\pmb{Aut}(D[p^l])^{G,0}_{\text{crys},\re}$ and $\Omega_l^{\mathfrak{g}}=\pmb{End}(D[p^l])^{\mathfrak{g},0}_{\text{crys},\re}$. Thus (a) follows from the property \ref{inv} (iv). 

From (a) and the property \ref{Lie} (ii) we get that the sequence $(\gamma_D^G(l))_{l\ge n_D+1}$ is constant. To prove (b) we can assume $n_D^G\ge 2$. From the last two sentences and Proposition \ref{P2} (d), we get that $n_D^G\le n_D+1$ (i.e. (b) holds). 
\endproof

\begin{thm}\label{T7}
We assume that $p$ is odd and that $\upsilon_D^G\le {{p-1}\over 2}$ (for instance, these hold if $2\le 2r<p$). If $n_D=0$, then we assume that $G_k$ is connected. Then the seven properties of Theorem \ref{T6} hold. 
\end{thm}

\noindent
{\bf Proof:} As in the proof of Theorem \ref{T6}, it suffices to show that the statement (a) of Theorem \ref{T6} holds. We take $\Sigma:=\Sigma_{\upsilon_D^G}$. Thus we have an isomorphism $\pmb{Aut}(D[p^m])^{G,0}_{\text{crys},\re}\to \pmb{End}(D[p^m])^{\mathfrak{g},0}_{\text{crys},\re}$, cf. property \ref{inv} (v). Therefore the statement (a) of Theorem \ref{T6} holds.
\endproof

\subsection{Remarks}\label{Remark2}

We end this section with three remarks.

\smallskip
{\bf (a)} Subsubsections \ref{tildeH} and \ref{digr} can be as well adapted to the relative context provided by $G$. For instance, the group scheme over $\Spec W(k)$ defined by the intersection $\pmb{Aut}(\phi,\vartheta)^G:=\pmb{Aut}(\phi,\vartheta)\cap (G^{(\sigma)}\times_{W(k)} G)$ taken inside $\GL_{M^{(\sigma)}}\times_{W(k)} \GL_M$, is naturally isomorphic to $(\tilde{\mathcal H}^G)^{(\sigma)}$. One ends up with an action $\Bbb T_m^{G,\prime:}\mathcal H_m^{G,\prime}\times_k \mathcal D_m^{G,\prime}\to \mathcal D_m^{G,\prime}$ isomorphic to the action $\Bbb T_m^{G,(\sigma)}$.

\smallskip
{\bf (b)} Referring to Properties \ref{inv} (iv) and (v), their isomorphisms $\Omega_m^G\to\Omega_m^{\mathfrak g}$ and $\pmb{Aut}(D[p^m])^{G,0}_{\text{crys},\re}\to \pmb{End}(D[p^m])^{\mathfrak{g},0}_{\text{crys},\re}$ are based on the closed subgroup scheme $G$ of $\GL_M$. One gets variants of these isomorphisms working with either the closed subgroup scheme $\phi^{-1}G\phi$ of $\GL_{\phi^{-1}(M)}$ or the closed subgroup scheme $\tilde{\mathcal H}^G$ of $\GL_{M,\phi^{-1}(M)}$ or the closed subgroup scheme $\pmb{Aut}(\phi,\vartheta)^G$ of $\GL_{M^{(\sigma)}}\times_{W(k)} \GL_M$. For instance, working with $\tilde{\mathcal H}^G$, the role played by $\pmb{Aut}(D[p^m])_{\text{crys},\re}$ is played by $\mathcal C^G_m$ and the role played by $\pmb{End}(D[p^m])_{\text{crys},\re}$ is played by the reduced closed subgroup scheme of $\Bbb W_m(\underline{\Lie(\tilde{\mathcal H}^G}))$ whose $k$-valued points are elements $e_m\in \Lie(\tilde{\mathcal H}^G)/p^m\Lie(\tilde{\mathcal H}^G)$ which have a lift $e\in  \Lie(\tilde{\mathcal H}^G)$ such that we have $\phi(e)-e\in p^m\mathfrak{g}$ (here we view $\Lie(\tilde{\mathcal H}^G)$ as a Lie subalgebra of $\mathfrak{g}$ via $\Lie(\tilde{\mathcal P}^G)$). Moreover, all these isomorphisms are compatible with the different isomorphisms and homomorphisms between the four groups $G$, $\phi^{-1} G\phi$, $\tilde{\mathcal H}^G$, $\pmb{Aut}(\phi,\vartheta)^G$, and their pull-backs via $\sigma$, as one can easily check based on Fact \ref{F5} of Appendix B. 

\smallskip
{\bf (c)} Axiom (AX2) was used only in the proofs of Lemmas \ref{L5} and \ref{L6} and Proposition \ref{P2} (a). Thus all the results of this section continue to hold even if the axiom (AX2) is replaced by the following weaker version of it.

\medskip
{\bf (AX$2'$)} The intersection $\mathcal{W}^G=G_k\cap \mathcal{W}$ is a smooth subgroup of $G_k$ and $G_{B(k)}$ is well positioned with respect to $(M,\phi^{-1}(M))$.

\section{Examples}\label{Sect7}

In this section we exemplify the results of Section \ref{Sect6} through three simple situations that are of special interest. Expressions of the form $\phi_1$ and $\vartheta_1$ will have a meaning different from the one before Subsection \ref{SSect3.1}.

\subsection{The homomorphism context}\label{E1}
Let $D_1$ and $D_2$ be two $p$-divisible groups over $k$. For $j\in\{1,2\}$ let $(M_j,\phi_j,\vartheta_j)$ be the contravariant Dieudonn\'e module of $D_j$. We take $D:=D_1\oplus D_2$. We have $(M,\phi,\vartheta)=(M_1,\phi_1,\vartheta_1)\oplus (M_2,\phi_2,\vartheta_2)$ and a direct sum decomposition $\End(M)=\End(M_1)\oplus \End(M_2)\oplus\Hom(M_1,M_2)\oplus \Hom(M_2,M_1)$
into $W(k)$-modules. Let $G$ be the smooth integral closed subgroup scheme of $\GL_M$ whose group of valued points in a commutative $W(k)$-algebra $R$ is 
$$G(R)=1_{R\otimes_{W(k)} M}+R\otimes_{W(k)} \Hom(M_2,M_1).$$ 
We have $\mathfrak{g}:=\Lie(G)=\Hom(M_2,M_1)$. It is easy to see that axioms (AX1) and (AX2) hold for the triple $(M,\phi,G)$ (in connection to (AX2) one has to take $M=F^1\oplus F^0$ such that we have $F^s=(F^s\cap M_1)\oplus (F^s\cap M_2)$ for each $s\in\{0,1\}$). Moreover, the condition  (i) of Theorem \ref{T6} holds. 

Let $n_{D_1,D_2}$ be the smallest non-negative integer for which the following statement holds. If $E$ is a $p$-divisible group over $k$ equipped with a short exact sequence $0\to D_2\to E\to D_1\to 0$ whose truncation of level $n_{D_1,D_2}$ splits, then the short exact sequence $0\to D_2\to E\to D_1\to 0$ splits. The contravariant Dieudonn\'e module of every $p$-divisible group $\tilde E$ over $k$ equipped with a short exact sequence $0\to D_2\to \tilde E\to D_1\to 0$ is isomorphic to $(M,g\phi,\vartheta g^{-1})$ for some arbitrary element $g\in G(W(k))$. Based on this one easily checks that $n_{D_1,D_2}=n_D^G$ (to be compared with [Va1, Lemma 3.2.2]). 

For $i\in\Bbb N$ let $\pmb{Hom}(D_1[p^i],D_2[p^i])$ be the group scheme over $k$ of homomorphisms from $D_1[p^i]$ to $D_2[p^i]$. Let $\gamma_{D_1,D_2}(i):=\dim(\pmb{Hom}(D_1[p^i],D_2[p^i]))$. One can naturally identify $\pmb{Hom}(D_1[p^i],D_2[p^i])_{\re}$ with $\pmb{End}(D)^{\mathfrak{g}}_{\re}$ as well as $\pmb{Aut}(D[p^i])^G_{\re}$ with $1_{D[p^i]}+\pmb{Hom}(D_1[p^i],D_2[p^i])_{\re}$. Therefore $\gamma_D^G(i)=\gamma_{D_1,D_2}(i)$. Thus from Theorem \ref{T6} and Proposition \ref{P2} (c) we get that:

\medskip
{\bf (i)} If $n_{D_1,D_2}>0$, then $\gamma_{D_1,D_2}(1)<\gamma_{D_1,D_2}(2)<\cdots<\gamma_{D_1,D_2}(n_{D_1,D_2})$. 

\smallskip
{\bf (ii)} The sequence $(\gamma_{D_1,D_2}(i+1)-\gamma_{D_1,D_2}(i))_{i\in \Bbb N}$ is a decreasing sequence in $\Bbb N$ and for $i\ge  n_{D_1,D_2}$ we have $\gamma_{D_1,D_2}(i)=\gamma_{D_1,D_2}(n_{D_1,D_2})$. Moreover we have $n_{D_1,D_2}\le n_{D_1\oplus D_2}$.

\smallskip
{\bf (iii)} For integers $l>i>0$, the image of the restriction homomorphism $\pmb{Hom}(D_1[p^l],D_2[p^l])\to \pmb{Hom}(D_1[p^i],D_2[p^i])$ is finite if and only if $l-i\ge n_{D_1,D_2}$. 

\subsection{The filtered context}\label{E2}
We assume that $D$ has an increasing filtration $\mathcal F$ of the form
$$0=D_0\subset D_1\subset\cdots \subset D_{s-1}\subset D_s=D.$$
Let $M_i$ be the direct summand of $M$ such that $M/M_i$ defines the contravariant Dieudonn\'e module of $D_i$. We have a decreasing filtration 
$$0=M_s\subset \cdots \subset M_1\subset M_0=M.$$
Let $G$ be the parabolic subgroup scheme of $\GL_M$ that normalizes each $M_j$ with $j\in\{0,\ldots,s\}$. We have $\mathfrak{g}:=\Lie(G)=\{e\in\End(M)|e(M_j)\subset M_j\;\forall j\in\{0,\ldots,s\}\}$. It is easy to see that axioms (AX1) and (AX2) hold for the triple $(M,\phi,G)$. Moreover, the condition (i) of Theorem \ref{T6} holds. 

Let $n_{\mathcal F}$ be the smallest non-negative integer for which the following statement holds. If $D^\prime$ is a $p$-divisible group over $k$ equipped with an increasing filtration $\mathcal F^\prime$ of the form
$0=D_0^\prime\subset D_1^\prime\subset\cdots \subset D_{s-1}^\prime\subset D_s^\prime=D^\prime$ whose truncation of level $n_{\mathcal F}$ is isomorphic to the truncation of level $n_{\mathcal F}$ of $\mathcal F$ and whose factors $D_{j+1}^\prime/D_j^\prime$ have the same dimension and codimension as $D_{j+1}/D_j$ for all $j\in\{0,\ldots,s-1\}$, then $\mathcal F^\prime$ is isomorphic to $\mathcal F$. The contravariant Dieudonn\'e module of each $p$-divisible group $D^\prime$ equipped with an increasing filtration $\mathcal F^\prime$ of the form
$0=D_0^\prime\subset D_1^\prime\subset\cdots \subset D_{s-1}^\prime\subset D_s^\prime=D^\prime$ whose factors $D_{j+1}^\prime/D_j^\prime$ have the same dimension and codimension as $D_{j+1}/D_j$ for all $j\in\{0,\ldots,s-1\}$, is isomorphic to $(M,g\phi,\vartheta g^{-1})$ for some arbitrary element $g\in G(W(k))$. Based on this one easily checks that $n_{\mathcal F}=n_D^G$ (to be compared with [Va1, Lemma 3.2.2]). 

For $i\in\Bbb N$, let $\pmb{Aut}_{\mathcal F}(D[p^i])$ (resp. $\pmb{End}_{\mathcal F}(D[p^i])$) be the group scheme of automorphisms (resp. endomorphisms) of the filtered truncated Barsotti--Tate group $0=D_0[p^i]\subset D_1[p^i]\subset\cdots \subset D_{s-1}[p^i]\subset D_s[p^i]=D[p^i]$.
Thus $\pmb{Aut}_{\mathcal F}(D[p^i])$ (resp. $\pmb{End}_{\mathcal F}(D[p^i])$) is a closed subgroup scheme of $\pmb{Aut}(D[p^i])$ (resp. $\pmb{End}(D[p^i])$) with the property that for a commutative $k$-algebra $R$ and for $e_R\in \pmb{Aut}(D[p^i])(R)$ (resp. $e_R\in \pmb{End}(D[p^i])(R)$) we have $e_R\in \pmb{Aut}_{\mathcal F}(D[p^i])(R)$ (resp. $e_R\in \pmb{End}_{\mathcal F}(D[p^i])(R)$) if and only if $e_R$ maps $D_j[p^i]_R$ into itself for all $j\in\{0,\ldots,s\}$. Each $\pmb{Aut}_{\mathcal F}(D[p^i])$ is an open (as well as closed) subscheme of $\pmb{End}_{\mathcal F}(D[p^i])$ and thus we can define $\gamma_{\mathcal F}(i):=\dim(\pmb{Aut}_{\mathcal F}(D[p^i]))=\dim(\pmb{End}_{\mathcal F}(D[p^i]))$. It is easy to see that $\pmb{Aut}_{\mathcal F}(D[p^i])_{\re}=\pmb{Aut}(D[p^i])^G_{\re}$ and $\pmb{End}_{\mathcal F}(D[p^i])_{\re}=\pmb{End}(D[p^i])^{\mathfrak{g}}_{\re}$. Thus $\gamma_{\mathcal F}(i)=\gamma_D^G(i)$. From Theorem \ref{T6} and Proposition \ref{P2} (c) we get that:

\medskip
{\bf (i)} If $n_{\mathcal F}>0$, then we have $\gamma_{\mathcal F}(1)<\gamma_{\mathcal F}(2)<\cdots<\gamma_{\mathcal F}(n_{\mathcal F})$. 

\smallskip
{\bf (ii)} The sequence $(\gamma_{\mathcal F}(i+1)-\gamma_{\mathcal F}(i))_{i\in \Bbb N}$ is a decreasing sequence in $\Bbb N$ and for $i\ge  n_{\mathcal F}$ we have $\gamma_{\mathcal F}(i)=\gamma_{\mathcal F}(n_{\mathcal F})$. Moreover we have $n_{\mathcal F}\le n_D$.

\smallskip
{\bf (iii)} For integers $l> i> 0$, the image of the restriction homomorphism $\pmb{End}_{\mathcal F}(D[p^l])\to \pmb{End}_{\mathcal F}(D[p^i])$ is finite if and only if $l-i\ge n_{\mathcal F}$. 

\smallskip
{\bf (iv)} If $i\ge 2n_{\mathcal F}$ is an integer, then the unipotent group $\pmb{Aut}_{\mathcal F}(D[p^i])^0_{\re}$ is commutative.

\subsection{The principally quasi-polarized context}\label{E3}

Let $D^{\text{t}}$ be the Cartier dual of $D$. We assume that there exists an isomorphism $\lambda:D\to D^{\text{t}}$ such that the perfect bilinear form $\psi:M\times M\to W(k)$ associated naturally to $\lambda$ is symmetric (resp. is alternating). If $p=2$ and $\psi$ is symmetric, we also assume that $\psi$ modulo $2$ is alternating. We refer to $\lambda$ as a symmetric principal (resp. a principal) quasi-polarization of $D$. 

We have $c=d>0$. Let $G:=\pmb{SO}(q)$, where $q:M\to W(k)$ is the quadratic form defined by the rule $q(x)={1\over 2}\psi(x,x)$ (resp. let $G:=\pmb{Sp}(M,\psi)$). The group scheme $G$ is $\Bbb G_m$ (resp. is $\pmb{SL}_2$) if $d=1$ and it is semisimple if $d>1$ (see [Va5, Subsect. 3.1 and Prop. 3.4] for the case when $p=2$ and $\psi$ is symmetric). It is easy to see that axioms (AX1) and (AX2) hold for the triple $(M,\phi,G)$. 

Let $n_{D,\lambda}$ be the smallest non-negative integer $n$ for which the following statement holds: if $(D^\prime,\lambda^\prime)$ is another $p$-divisible group of dimension $d$ over $k$ equipped with a symmetric principal (resp. with a principal) quasi-polarization such that $(D^\prime[p^n],\lambda^\prime[p^n])$ is inner isomorphic (resp. is isomorphic) to $(D[p^n],\lambda[p^n])$, then $(D^\prime,\lambda^\prime)$ is isomorphic to $(D,\lambda)$. Here, in the case when $\psi$ is symmetric, by an inner isomorphism $(D^\prime[p^n],\lambda^\prime[p^n])\to (D[p^n],\lambda[p^n])$ we mean an isomorphism defined at the level of Dieudonn\'e modules by an isomorphism $(M/p^nM,\phi_n,\vartheta_n,q_n)\to (M^\prime/p^nM^\prime,\phi_n^\prime,\vartheta_n^\prime,q^\prime_n)$, where the left lower indexes $n$ mean reduction modulo $p^n$. 

We always have $n_{D,\lambda}\ge 1$. We have $n_D^G=0$ if and only if $d=1$ and $\psi$ is symmetric. If $n_D^G\ge 1$, it is well known that $n_{D,\lambda}=n_D^G$ (if $\psi$ is alternating, see [Va2, Example 4.5]; the symmetric case is similar). For $i\in\Bbb N$, let $\gamma_{D,\lambda}(i)$ be the dimension of the group scheme $\pmb{Inn}(D[p^i],\lambda[p^i])$ (resp. $\pmb{Aut}(D[p^i],\lambda[p^i])$) of inner automorphisms (resp. of automorphisms) of $(D[p^i],\lambda[p^i])$, where $\pmb{Inn}(D[p^i],\lambda[p^i])$ is  some  subgroup scheme of $\pmb{Aut}(D[p^i])$ whose $k$-valued points are the inner automorphisms defined above. It is easy to see that $\pmb{Inn}(D[p^i],\lambda[p^i])_{\re}=\pmb{Aut}(D[p^i])^G_{\re}$ (resp. $\pmb{Aut}(D[p^i],\lambda[p^i])_{\re}=\pmb{Aut}(D[p^i])^G_{\re}$). Thus $\gamma_{D,\lambda}(i)=\gamma_D^G(i)$. 

If $p>2$, then $G$ admits a Cayley transform of the form $g\mapsto (1_M-g)(1_M+g)^{-1}$ and thus Theorem \ref{T6} applies. In particular, for $p>2$ and $d>1$ (resp. for $p>2$) we have $n_{D,\lambda}=n_D^G\le n_D$. If $p=2$ and $d>1$ (resp. if $p=2$), then from Proposition \ref{P3} (b) we get the weaker fact that $n_{D,\lambda}=n_D^G\le n_D+1$. 

If $p>2$ and $d>1$ (resp. if $p>2$), then from Theorem \ref{T6} we also get that:

\medskip
{\bf (i)} We have $\gamma_{D,\lambda}(1)<\gamma_{D,\lambda}(2)<\cdots<\gamma_{D,\lambda}(n_{D,\lambda})$. 

\smallskip
{\bf (ii)} The sequence $(\gamma_{D,\lambda}(i+1)-\gamma_{D,\lambda}(i))_{i\in \Bbb N}$ is a decreasing sequence in $\Bbb N$ and for $i\ge  n_{D,\lambda}$ we have $\gamma_{D,\lambda}(i)=\gamma_{D,\lambda}(n_{D,\lambda})$. 

\smallskip
{\bf (iii)} For integers $l>i>0$, the restriction homomorphism $\pmb{Inn}(D[p^l],\lambda[p^l])\break\to \pmb{Inn}(D[p^i],\lambda[p^i])$ (resp. $\pmb{Aut}(D[p^l],\lambda[p^l])\to \pmb{Aut}(D[p^i],\lambda[p^i])$) has a finite image if and only if $l-i\ge n_{D,\lambda}$. 

\smallskip
{\bf (iv)} If $i\ge 2n_{D,\lambda}$ is an integer, then the unipotent group $\pmb{Inn}(D[p^i],\lambda[p^i])_{\re}^0$ (resp. $\pmb{Aut}(D[p^i],\lambda[p^i])_{\re}^0$) is commutative.  

\bigskip
The notations of the appendices below will be independent from the previous notations of the paper. 

\section{Appendix A: affine group schemes over discrete valuation rings}

Let $R$ be an arbitrary discrete valuation ring. Let $K$ be the field of fractions of $R$ and let $\pi$ be a uniformizer of $R$. Let $k:=R/(\pi)$ be the residue field of $R$. Let $H=\Spec Q$ be a flat, affine group scheme over $\Spec R$. Let $J=\Spec Q/I$ be a closed subgroup scheme of $H_k$; thus $I$ is an ideal of $Q$ which contains $Q\pi$. Let $Q'$ be the $R$-subalgebra of $Q[{1\over {\pi}}]$ generated by all elements of the form ${x\over {\pi}}$ with $x\in I$; it contains $Q$. See [BLR, Ch. 3, Sect. 3.2] for dilatations. The dilatation of $H$ centered on $J$ is the flat, affine group scheme $H'=\Spec Q'$ over $\Spec R$. One has the following two properties (they follow directly from the very description of $Q'$; see [BLR, Ch. 3, Sect. 3.2, Props. 1, 2, and 3]):

\medskip
{\bf (i)} there exists a homomorphism $f:H'\to H$ that corresponds to the $R$-monomorphism $Q\hookrightarrow Q'$; 

\smallskip
{\bf (ii)} a morphism $g:X\to H$ of flat $\Spec R$-schemes factors
 through $H'$ if and only if the morphism $g_k:X_k\to H_k$ factors through $J$, and the factorization is unique if it exists.

\medskip
By applying (ii) with $X$ a $\Spec K$-scheme, we get that $f_K$ is an isomorphism of $\Spec K$-schemes that corresponds to the identity $Q[{1\over {\pi}}]=Q'[{1\over {\pi}}]$. If $H$ and $J$ are smooth, then $H'$ is also smooth (cf. loc. cit.). If $\tilde H$ is a flat, closed subgroup scheme of $H$ and if $\tilde J:=J\cap \tilde H_k$, then the dilatation $\tilde H'$ of $\tilde H$ centered on $\tilde J$ comes equipped with a homomorphism $\tilde H'\to H'$ (cf. (ii)). It is easy to see (cf. loc. cit.) that:

\medskip
{\bf (iii)} the homomorphism $\tilde H'\to H'$ is a closed embedding.

\medskip
Let $V$ be a finite dimensional $K$-vector space. By an $R$-{\it lattice} of $V$ we mean an $R$-submodule of $V$ generated by a $K$-basis of $V$. Let $M$ and $M'$ be two lattices of $V$ such that we have $\pi M'\subset M\subset M'$. We consider a direct sum decomposition $M=F^1\oplus F^0$ into free $R$-modules such that we have $M'=\pi^{-1}F^1\oplus F^0$. Let $\lambda:\Bbb G_m\to\GL_V$ be the cocharacter defined by the rule: $t\in \Bbb G_m(K)$ fixes $F^0[{1\over {\pi}}]$ and acts as the multiplication by $t^{-1}$ on $F^1[{1\over {\pi}}]$.

Let $G$ be a closed subgroup scheme of $\GL_V$. Let $G_M$ and $G_{M'}$ be the schematic closures of $G$ in $\GL_M$ and $\GL_{M'}$ (respectively); they are flat, affine group schemes over $\Spec R$. 
Let $G_{M,M'}$ be the schematic closure of $G$ in $\GL_M\times_R \GL_{M'}$ embedded diagonally into the generic fibre of $\GL_M\times_R \GL_{M'}$. 

\medskip\noindent
{\bf Example} Suppose $G=\GL_V$. Then $G_M=\GL_M$, $G_{M'}=\GL_{M'}$, and we claim that $G_{M,M'}$ is the group scheme $H$ over $\Spec R$ that represents the functor of automorphisms of the pair of $R$-linear maps $\xymatrixcolsep{0.8pc}\xymatrix@1{ M \ar@<0.5ex>[rr]^{\iota} & &M'
\ar@<0.5ex>[ll]^{\pi} \\}$, where $\iota:M\hookrightarrow M'$ is the inclusion. As $\lambda(\pi^{-1})$ maps $M'$ isomorphically onto $M$, $H$ is isomorphic to the group scheme $\pmb{Aut}(\iota_1,\iota_2)$ over $\Spec R$ that represents the functor of automorphisms of the pair of $R$-linear maps $\xymatrixcolsep{0.8pc}\xymatrix@1{ M \ar@<0.5ex>[rr]^{\iota_1} & &M \ar@<0.5ex>[ll]^{\iota_2} \\}$, where $\iota_1:M\hookrightarrow M$ and $\iota_2:M\hookrightarrow M$ map $(x,y)\in F^1\oplus F^0=M$ to $(\pi x,y)$ and $(x,\pi y)$ (respectively). For a commutative $R$-algebra $S$, we write an $S$-linear
endomorphism of $S\otimes_R (F^1\oplus F^0)$ in the form $\binom{s\; t}{u\; v}$,
where $s\in \End(S\otimes_R F^1)$, $t\in \Hom(S\otimes_R F^0,S\otimes_R
F^1)$, $u\in \Hom(S\otimes_R F^1,S\otimes_R F^0)$, and $v\in
\End(S\otimes_R F^0)$. The subgroup $\pmb{Aut}(\iota_1,\iota_2)(S)$ of
$\GL_{M}(S)\times\GL_M(S)$ is formed by all pairs of the form
$(\binom{s\;\;\; t}{\pi u\; v},\binom{s\; \pi t}{u\;\; v})$ which are invertible. From this it follows that $H$ is flat (in fact even smooth) and thus it is $G_{M,M'}$. The description of the subgroups $\pmb{Aut}(\iota_1,\iota_2)(S)$ also implies that the following three properties hold: 

\medskip
{\bf (a)} the first projection homomorphism $\pmb{Aut}(\iota_1,\iota_2)\to \GL_M$ and thus also $\GL_{M,M'}\to \GL_M$ is the dilatation of $\GL_M$ centered on the parabolic subgroup $P_k$ of $\GL_{M/\pi M}$ that normalizes the $k$-vector subspace $\pi M'/\pi M=F^1/\pi F^1$ of $M/\pi M$;

\smallskip
{\bf (b)}  similarly, the second projection homomorphism $\pmb{Aut}(\iota_1,\iota_2)\to \GL_M$ (resp. $\GL_{M,M'}\to \GL_{M'}$) is the dilatation of $\GL_M$ (resp. of $\GL_{M'}$) centered on the parabolic subgroup $P'_k$ of $\GL_{M/\pi M}$ (resp. of $\GL_{M'/\pi M'}$) that normalizes the $k$-vector subspace $F^0/\pi F^0$ of $M/\pi M$ (resp. of $M'/\pi M'$);

\smallskip
{\bf (c)} the special fibre $\pmb{Aut}(\iota_1,\iota_2)_k$ (and thus also $H_k$) is connected. 

\medskip
From this Example and the property (iii) we get directly that:

\begin{fact}\label{F3}
In general, the homomorphism $G_{M,M'}\to G_M$ is the dilatation of $G_M$ centered on $G_{M,k}\cap P_k$ and the homomorphism $G_{M,M'}\to G_{M'}$ is the dilatation of $G_{M'}$ centered on $G_{M',k}\cap P'_k$.
\end{fact}  

\medskip
We also consider the intersection $G^*_{M,M'}:=(G_M\times_R G_{M'})\cap \GL_{M,M'}$ taken inside $\GL_M\times_R \GL_{M'}$. We have a closed embedding homomorphism $h:G_{M,M'}\to G^*_{M,M'}$ whose generic fibre is an isomorphism.

\begin{df}\label{df3}
We say that $G$ is well positioned with respect to $(M,M')$ if $G_{M,M'}=G^*_{M,M'}$, i.e. $h$ is an isomorphism (equivalently, $G^*_{M,M'}$ is flat).
\end{df}

For $n\in\Bbb N^*$ let $h_n:=h(R/(\pi^n)):G_{M,M'}(R/(\pi^n))\to G^*_{M,M'}(R/(\pi^n))$. If $G$ is well positioned with respect to $(M,M')$, then each $h_n$ is an isomorphism. 

\begin{fact}\label{F4}
We assume that $k$ is algebraically closed and $G_{M,M'}$ is smooth. If $h_1$ and $h_2$ are isomorphisms, then $G$ is well positioned with respect to $(M,M')$.
\end{fact}

\noindent
{\bf Proof:} By the fiberwise criterion of flatness, it suffices to check that the special fibre $h_k:G_{M,M',k}\to G^*_{M,M',k}$ of $h$ is an isomorphism. Since $h_k$ is a bijection (as $h_1$ is an isomorphism) and it induces an isomorphism $dh_k$ at the level of Lie algebras (as $h_2$ is an isomorphism), it is a closed embedding with smooth source and target. Thus, $h_k$ is an isomorphism.\endproof

\begin{prop}\label{P4}
If $G$ is normalized by $\lambda$, then $G$ is well positioned with respect to $(M,M')$.
\end{prop}

\noindent
{\bf Proof:} We will use the theory of dynamically defined (see [CGP, Sect. 2.1]) closed subgroup schemes $U_{*}(-\lambda)$, $Z_{*}(\lambda)$, and $U_{*}(\lambda)$ of $*$, where $*$ stands for an affine group scheme over either $\Spec K$ or $\Spec R$ normalized by $\lambda$ or by some extension of $\lambda$ to a cocharacter over $\Spec R$. For example, $Z_{*}(\lambda)$ is the scheme theoretic centralizer of $\lambda$ (or of some extension of it over $\Spec R$) in $*$, while $U_{*}(\lambda)$ and $U_{*}(-\lambda)$ are unipotent subgroup schemes of $*$, which for $*=\GL_V$ are obtained naturally from the $K$-vector spaces $\Hom(F^1[{1\over {\pi}}],F^0[{1\over {\pi}}])$ and $\Hom(F^0[{1\over {\pi}}],F^1[{1\over {\pi}}])$ (respectively). 

The product open embedding $U_{\GL_V}(-\lambda)\times_K Z_{\GL_V}(\lambda)\times_K U_{\GL_V}(\lambda)\hookrightarrow \GL_V$ defines the ``big cell'' of $\GL_V$. Its intersection with $G$ is the ``big cell'' $\mathfrak{C}_G:=U_G(-\lambda)\times_K Z_G(\lambda)\times_K U_G(\lambda)$ of $G$, cf. [CGP, Prop. 2.1.8 and Rm. 2.1.11]. The ``big cell'' notion extends to smooth affine group schemes over $R$ (cf. loc. cit.) and the special fibre of $\GL_{M,M'}$ lies in the big cell. When we consider the big cell $\mathfrak{C}_M\times_R \mathfrak{C}_{M'}$ of $\GL_{M}\times_R \GL_{M'}$, it is easy to see that we have corresponding decompositions for $G_M$, $G_{M'}$, and $G_{M,M'}$ and thus also for $G^*_{M,M'}=(G_M\times_R G_{M'})\cap \GL_{M,M'}$. For instance for $G_M$, if $\mathfrak{C}_{G_M}$ is the schematic closure of $\mathfrak{C}_G$ in $\mathfrak{C}_M$, then the closed embedding morphism $j:\mathfrak{C}_{G_M}\hookrightarrow G_M\cap \mathfrak{C}_M$ between flat $\Spec R$-schemes is an isomorphism as $j_K$ is an isomorphism; but $\mathfrak{C}_{G_M}$ is a product $U_{G_M}(-\lambda)\times_R Z_{G_M}(\lambda)\times_R U_{G_M}(\lambda)$ (as its generic fibre $\mathfrak{C}_G$ and $\mathfrak{C}_M$ are such products) and thus we get the product decomposition $G_M\cap \mathfrak{C}_M=U_{G_M}(-\lambda)\times_R Z_{G_M}(\lambda)\times_R U_{G_M}(\lambda)$ into $\Spec R$-schemes which are flat as $G_M$ is so. Thus it suffices to check the flatness of $G^*_{M,M'}$ for each one of the three factors of its open subscheme big cell. 

{\bf The $Z(\lambda)$ factor.} We have $Z_G(\lambda)\subset \GL_{F^0[{1\over {\pi}}]}\times_K \GL_{F^1[{1\over {\pi}}]}$. Its schematic closures in $G_M$ and $G_{M'}$ are $Z_{G_M}(\lambda)$ and $Z_{G_{M'}}(\lambda)$ (respectively) and can be identified with the schematic closures of $Z_G(\lambda)$ in $Z_{\GL_M}(\lambda)=\GL_{F^1}\times_R \GL_{F^0}$ and $Z_{\GL_{M'}}(\lambda)=\GL_{\pi^{-1}F^1}\times_R \GL_{F^0}$ (respectively). As $\GL_{F^1}=\GL_{\pi^{-1}F^1}$, $Z_{G^*_{M,M'}}(\lambda)$ is the intersection of the product of these last schematic closures with the diagonal of $Z_{\GL_M}(\lambda)\times_R Z_{\GL_{M'}}(\lambda)$ and this intersection is flat and can be identified  with both $Z_{G_M}(\lambda)$ and $Z_{G_{M'}}(\lambda)$.

{\bf The $U(\lambda)$ factor.} We have $U_{\GL_M}(\lambda)=\underline{Hom(F^1,F^0)}$ and $U_{\GL_{M'}}(\lambda)=\underline{Hom(\pi^{-1}F^1,F^0)}$ (here underlying means the group scheme functor associated to a module). Moreover $U_{\GL_{M,M'}}(\lambda)$ is the graph of the restriction homomorphism $\phi:\underline{Hom(\pi^{-1}F^1,F^0)}\to \underline{Hom(F^1,F^0)}$. But $\phi$ sends the schematic closure $U_{G_{M'}}(\lambda)$ of $U_G(\lambda)$ in $\underline{Hom(\pi^{-1}F^1,F^0)}$ to the schematic closure  $U_{G_M}(\lambda)$ of $U_G(\lambda)$ in $\underline{Hom(F^1,F^0)}$. Thus again the intersection $U_{G^*_{M,M'}}(\lambda)$ of the product of these schematic closures with the graph of $\phi$ is flat and in fact it is isomorphic to $U_{G_{M'}}(\lambda)$. 

{\bf The $U(-\lambda)$ factor.} The argument for this is similar to the one in the previous paragraph and in fact we have $U_{G^*_{M,M'}}(-\lambda)=U_{G_M}(-\lambda)$.\endproof

\section{Appendix B: properties of $\exp$ and $\log$ maps for matrices in $p$-adic contexts}

Let $p$ be a prime. Let $K$ be a field of characteristic $0$ complete with respect to an ultrametric absolute value $|\;|:K\to\Bbb R$. We extend $|\;|$ to an absolute value on an algebraic closure $\bar K$ of $K$, denoted also by $|\;|$. We assume that the residue characteristic of $K$ is $p$ and let $q:=|p|^{1\over {p-1}}\in (0,1)$. Let $n\in \Bbb N^{\ast}$ and let $I_n$ be the identity matrix of $M_n(K)$ or $M_n(\bar K)$. By the norm of a matrix $X=(x_{i,j})_{1\le i,j\le n}\in M_n(\bar K)$ we mean
$$||X||:=\max\{|x_{i,j}||1\le i,j\le n\}\in [0,\infty).$$
Let $\rho(X)$ (denoted $||X||_{\text{sp}}$ in some papers) be the {\it spectral radius} of $X$, i.e. the maximum of the absolute values of the eigenvalues of $X$ in $\bar K$. We can conjugate $X$ to an upper triangular matrix whose entries off the main diagonal have arbitrarily small absolute values. Thus if $\rho(X)>0$, then we can conjugate $X$ so that we have $||X||=\rho(X)$ and if $\rho(X)=0$, then we can conjugate $X$ so that its norm is arbitrarily small. The series
$e_1(X)=\sum_{n=1}^{\infty} {{X^n}\over {n!}}$ and $\exp(X)=I_n+e_1(X)$ converge if $\rho(X)<q$ and the series $l_1(X)=\log(I_n+X)=\sum_{n=1}^{\infty} {{(-1)^{n-1}}\over n}X^n$ converges if $\rho(X)<1$. From [Bo, Ch. II, Sect. 8, Prop. 4] and our conjugation statements, we get that $e_1$ and $l_1$ define inverse bijections on $\Bbb B_n:=\{X\in M_n(K)|\rho(X)<q\}$. 

\begin{prop}\label{P5}
Let $G$ be an algebraic subgroup of $\GL_{n,K}$. Then we have inverse bijections $\xymatrixcolsep{0.9pc}\xymatrix@1{ \{\zeta\in\Lie(G)|\rho(\zeta)<q\} \ar@<0.5ex>[rr]^{\exp} & &\{g\in G(K)|\rho(g-I_n)<q\}\ar@<0.5ex>[ll]^{\log} \\}$.
\end{prop}
\noindent
{\bf Proof:} For $\zeta\in M_n(K)$ with $\rho(\zeta)<q$ we have to check that $\zeta\in\Lie(G)$ if and only if $\exp(\zeta)\in G(K)$.

To check the only if part, we recall that we have a formal exponential $\exp(T\zeta)$ in $G(K[[T]])$ (cf. [DGa, Ch. II, Sect. 6, Subsect. 3]). By Subsubsection 3.3 of loc. cit., it is given by the usual exponential series for $\GL_{n,K}$, hence the spectral radius condition gives that $\exp(T\zeta)$ and $\exp(-T\zeta)$ are matrices with coefficients in the ring  $K\{T\}$ of restricted formal power series. Thus  $\exp(T\zeta)\in G(K[[T]])\cap\GL_n(K\{T\})=G(K\{T\})$ and by evaluating at $T=1$ we get that $\exp(\zeta)\in G(K)$.

We now check the if part. We have $\exp(T\zeta)\in\GL_n(K\{T\})$ and it belongs to $G(K)$ when evaluated at $T=1$ and thus at each $T\in\Bbb Z$. Hence every equation for $G$ in  $\GL_{n,K}$ evaluated at $\exp(T\zeta)$ is an element in $K\{T\}$ with infinitely many zeros in $K$ and thus it is $0$, cf. Weierstra{\ss} preparation theorem for $K\{T\}$. Thus $\exp(T\zeta)\in G(K\{T\})\subset G(K[[T]])$ and by evaluating its derivative at $T=0$ we get that $\zeta\in\Lie(G)$.\endproof

\medskip\noindent
{\bf Remark} Proposition \ref{L6} and its proof extend to the case when the residue field of $K$ has characteristic $0$ provided we take $q=1$.

\subsection{Applications to rings of Witt vectors} 

Let $A$ be a reduced $\Bbb F_p$-algebra. Let $W(A)$ be the ring of ($p$-typical) Witt vectors with coefficients in $A$. Let $F$ and $V$ be the Frobenius and the Verschiebung maps of $W(A)$. Let $m\in\Bbb N^\ast$. Let $W_m(A):=W(A)/V^m(W(A))$. For $X=(x_{i,j})_{1\le i,j\le n}\in M_n(W(A))$ and $l\in\Bbb N$, let $\pi_l(X)=(x_{l,i,j})_{1\le i,j\le n}\in M_n(A)$ be the matrix of the $l$-th Witt components of the entries of $X$; thus $x_{i,j}=(x_{0,i,j},x_{1,i,j},\ldots)\in W(A)$. For $s\in\Bbb N^{\ast}$, let 
$$\Sigma_s=\{X\in M_n(W(A))|\pi_0(X)^s=0\}\;\;\text{and}\;\; \Sigma_{1,s}=V(\Sigma_s).$$
For $a\in\Bbb N^{\ast}$ we have a sequence of inclusions
$$\Sigma_{1,1}\subset\Sigma_{1,2}\subset\cdots\subset \Sigma_{1,n}=\Sigma_{1,n+a}\subset \Sigma_1\subset\Sigma_2\subset\cdots \subset\Sigma_n=\Sigma_{n+a}.$$ 
Moreover, $\Sigma_1=V(M_n(W(A)))$ and $\Sigma_{1,1}=V^2(M_n(W(A)))$.

As $A$ is reduced, $W(A)$ is torsion free and thus expressions of the form ${y\over d}$ are uniquely determined if they make sense (here $y\in W(A)$ and $d\in\Bbb N^{\ast})$. The $V$-adic filtration is compatible with the multiplication, i.e. we have $V^a(W(A))V^b(W(A))\subset V^{a+b}(W(A))$ for all $a,b\in\Bbb N$. Also we recall that $V(x_1)\cdots V(x_d)=p^{d-1}V(x_1\cdots x_d)$ for all $d\in\Bbb N^{\ast}$ and $x_1,\ldots,x_d\in W(A)$.  Moreover, as $V$ is additive, these properties extend to matrices in $M_n(W(A))$.

\begin{prop}\label{P6}
Let $\Sigma$ be either $\Sigma_s$ for $p>2$ and $1\le s\le {{p-1}\over 2}$ or $\Sigma_{1,t}$ for $p$ arbitrary and $t\in\Bbb N^*$. Then the following three properties hold:

\medskip
{\bf (a)} For $v\in\Bbb N^{\ast}$ the rule $X\mapsto {{X^v}\over {v!}}$ sends $\Sigma$ to $\Sigma$.

\smallskip
{\bf (b)} The partial sums of the series $e_1(\flat)=\sum_{v=1}^{\infty} {{\flat^v}\over {v!}}$ and $l_1(\flat)=\sum_{v=1}^{\infty} {{(-1)^{v-1}}\over v}\flat^v$ send $\Sigma$ to $\Sigma$ and converge $V$-adically to maps $e_1:\Sigma\to\Sigma$ and $l_1:\Sigma\to\Sigma$ which are inverse to each other.

\smallskip
{\bf (c)} For $l\in\Bbb N$ the entries of $\pi_l(e_1(X))$ and $\pi_l(l_1(X))$ are polynomials with coefficients in $\Bbb F_p$ in the entries of $\pi_{\tilde l}(X)$ with $\tilde l\in\{0,\ldots,l\}$. Thus if $\Sigma[m]:=\Image(\Sigma\to M_n(W_m(A)))$, then $e_1$ and $l_1$ induce inverse bijections $\xymatrixcolsep{0.8pc}\xymatrix@1{ \Sigma[m] \ar@<0.5ex>[rr]^{e_1} & &\Sigma[m]
\ar@<0.5ex>[ll]^{l_1} \\}$ which commute with the conjugation by $\pmb{GL}_n(W_m(A))$.
\end{prop}

\noindent
{\bf Proof:} We have ${{V(X)^v}\over {v!}}={{p^{v-1}}\over {v!}} V(X^v)$ and $\text{ord}_p(v!)\le {{v-1}\over {p-1}}$ and thus the rule $X\mapsto {{X^v}\over {v!}}$ sends $\Sigma_1$ to $\Sigma_1$ and $\Sigma_{1,t}$ is preserved. If $v\ge td$ (resp. $v\ge td+1$), then the rule $X\mapsto {{X^v}\over {v!}}$ sends $\Sigma_{1,t}$ to $V^d(\Sigma_1)$ (resp. to  $V^d(\Sigma_{1,t})$).

If $p>2$ and $X\in\Sigma_s$ with $s\le {{p-1}\over 2}$, we can write $X^{{p-1}\over 2}=V(Y)$ and thus ${{X^{p-1}}\over p}=V(Y^2)\in\Sigma_s$; as ${{T^v}\over {v!}}$ is $cT^a({{T^{p-1}}\over p})^b$ for some $c\in \Bbb Z_{(p)}$, $a\in\{0,\ldots,p-1\}$, and $b\in\Bbb N$, we get that ${{X^v}\over {v!}}\in\Sigma_s$. If $v\ge (p-1)(d+1)$, then the rule $X\mapsto {{X^v}\over {v!}}$ sends $\Sigma_{s}$ to $V^d(\Sigma_1)$. Thus (a) and the first part of (b) hold.

The fact that $e_1$ is inverse to $l_1$ can be reduced to the case when $A$ is a perfect field $k$ in which case Proposition \ref{P5} applies with $K=W(k)[{1\over p}]$, $\Sigma_s\subset \{X\in M_n(K)=\Lie(\GL_{n,K})|\rho(X)\le |p|^{1\over s}\}\subset\Bbb B_n$, and  $\Sigma_{1,t}\subset  \{X\in M_n(K)=\Lie(\GL_{n,K})|\rho(X)\le |p|^{{t+1}\over t}\}\subset \Bbb B_n$. Thus the second part of (b) also holds. Part (c) follows directly from the above facts on $X\mapsto {{X^v}\over {v!}}$.\endproof

\medskip
Let $P$ be a finitely generated, projective $W_m(A)$-module. If $P$ is free of rank $n$, we can speak about $\Sigma_P\subset \End(P)$ and about inverse bijections $\xymatrixcolsep{0.8pc}\xymatrix@1{ \Sigma_P \ar@<0.5ex>[rr]^{e_1} & &\Sigma_P
\ar@<0.5ex>[ll]^{l_1} \\}$ which can be identified with $\xymatrixcolsep{0.8pc}\xymatrix@1{ \Sigma[m] \ar@<0.5ex>[rr]^{e_1} & &\Sigma[m]\ar@<0.5ex>[ll]^{l_1} \\}$ with respect to each isomorphism $P\to W_m(A)^n$ of $W_m(A)$-modules. These inverse bijections make sense for all $P$ (without assuming their freeness) via a natural sheafification process (there exists a finite open cover $\Spec A=\cup_{i=1}^u \Spec A_{f_i}$ such that $W_m(A_{f_i})\otimes_{W_m(A)} P$ is free for all $i\in\{1,\ldots,u\}$; here each $f_i\in A$). 

\begin{fact}\label{F5} 
{\bf (a)} If $P$ is a finitely generated, projective $W_m(A)$-module and $Q$ is a direct summand of $P$ and if $x\in\Sigma_P$ is such that $x(Q)\subset Q$, then we have $e_1(x)(Q)\subset Q$ and $l_1(x)(Q)\subset Q$.

\smallskip
{\bf (b)} If $\phi:P_1\to P_2$ is a $W_m(A)$-linear map between finitely generated, projective $W_m(A)$-modules and if $x_1\in \Sigma_{P_1}$ and $x_2\in\Sigma_{P_2}$ are such that $x_2\circ\phi=\phi\circ x_1$, then we have $e_1(x_2) \circ \phi=\phi\circ e_1(x_1)$ and  $l_1(x_2) \circ \phi=\phi\circ l_1(x_1)$.
\end{fact}

\noindent
{\bf Proof:} Locally in the Zariski topology of $\Spec A$, to prove (a) we can assume that $P$ is free and $Q$ is generated by a subset of a $W_m(A)$-basis for $P$. It is clear from the definitions that $e_1$ and $l_1$ preserve block upper triangular matrices of the form  $\binom{\bigstar\; \bigstar}{0\; \bigstar}$ and therefore (a) holds. Part (b) follows from (a) applied to the direct summand of $P_1\oplus P_2$ which is the graph of $\phi$.\endproof

\subsection{The perfect field case} Let now $A$ be a perfect field $k$ of characteristic $p$. Let $B(k):=W(k)[{1\over p}]$. Let $M$ be a free $W(k)$-module of finite rank. For $x\in\End(M)$, let $\bar x\in\End(M/pM)$ be its reduction modulo $p$. Let $G$ be a flat, closed subgroup scheme of $\GL_M$. Let $\mathfrak{g}:=\Lie(G_{B(k)})\cap \End(M)$; thus $\mathfrak{g}/p^m\mathfrak{g}$ embeds in $\End(M/p^mM)$ for every $m\in\Bbb N^*$. We define $\Sigma$ subsets $\Sigma(M)$ of $\End(M)$ via the rules: if $\Sigma=\Sigma_s$, then $\Sigma(M):=\{x\in\End(M)|\bar x^s=0\}$ and if $\Sigma=\Sigma_{1,t}$, then $\Sigma(M):=p\{x\in\End(M)|\bar x^t=0\}$. 

 \begin{prop}\label{P7}
{\bf (a)} For every $\Sigma$ as in Proposition \ref{P6}, the $\exp$ and $\log$ induce inverse bijections
$\mathfrak{g}\cap \Sigma(M)\rightleftarrows \{g\in G(W(k))|g-1_M\in\Sigma(M)\}$
which in turn induce for each $m\in\Bbb N^{\ast}$ inverse bijections between the image $\Sigma^{\mathfrak{g}}[m]$ of $\mathfrak{g}\cap \Sigma(M)$ in $\mathfrak{g}/p^m\mathfrak{g}$ and the image $\Sigma^G[m]$ of $\{g\in G(W(k))|g-1_M\in\Sigma(M)\}$ in $G(W_m(k))$. 

{\bf (b)} If $(M,\phi,\vartheta)$ is the Dieudonn\'e module of a $p$-divisible group $D$ over $k$ (thus $\phi:W(k)\otimes_{F,W(k)} M\to M$ and $\vartheta:M\to W(k) \otimes_{F,W(k)} M$ are $W(k)$-linear maps such that $\phi\circ\vartheta$ and $\vartheta\circ\phi$ are the multiplication by $p$) and if $e_m\in \Sigma^{\mathfrak{g}}[m]$ is the crystalline realization of an endomorphism of $D[p^m]$ (i.e. is an endomorphism of $(M,\phi,\vartheta)$ modulo $p^m$), then $\exp(e_m)$ and $\log(1_{M/p^mM}+e_m)$ are crystalline realizations of endomorphisms of $D[p^m]$.
\end{prop}

\noindent
{\bf Proof:} As $V(W(k))=pW(k)$, part (a) follows from  Propositions \ref{P5} and \ref{P6}. Part (b) follows from Fact \ref{F5} (b) applied to the reductions modulo $p^m$ of $\phi$ and $\vartheta$.\endproof 

\smallskip\noindent
{\bf Acknowledgments.} The second author would like to thank Binghamton University, IHES, Bures-sur-Yvette, and IAS, Princeton
for good work conditions and J. E. Humphreys for the
reference to [Sp] used in the proof of Proposition \ref{P1}. He was partially supported by the NSF grant DMS \#0900967. The authors would like to thank the referee for many valuable comments and suggestions.

\hbox{}
\hbox{Ofer Gabber,\;\;\;Email: gabber@ihes.fr}
\hbox{Address: IH\'ES, Le Bois-Marie, 35, Route de Chartres,} 
\hbox{F-91440 Bures-sur-Yvette, France.}

\bigskip

\hbox{Adrian Vasiu,\;\;\;Email: adrian@math.binghamton.edu}
\hbox{Address: Department of Mathematical Sciences, Binghamton University,}
\hbox{Binghamton, P. O. Box 6000, New York 13902-6000, U.S.A.}
\end{document}